\documentclass[11pt,english]{article}
\usepackage[T1]{fontenc}
\usepackage[latin1]{inputenc}
\usepackage{amsthm}
\usepackage{amsmath}
\usepackage{amssymb}
\usepackage{amsfonts}
\usepackage[all]{xy}
\usepackage{babel}
\usepackage{hyperref}


\DeclareMathOperator{\Alb}{Alb}

\DeclareMathOperator{\CH}{CH}

\DeclareMathOperator{\Decf}{\underline{Dec}}

\DeclareMathOperator{\Divf}{\underline{Div}}

\DeclareMathOperator{\FDivf}{\underline{FDiv}}
\DeclareMathOperator{\LDiv}{LDiv}

\DeclareMathOperator{\Ext}{Ext}

\DeclareMathOperator{\Extabk}{Ext_{\Abk}}

\DeclareMathOperator{\Gal}{Gal}

\DeclareMathOperator{\Homk}{Hom_{\fld}}
\DeclareMathOperator{\Homck}{Hom_{\fld}^{\mathrm{cont}}}

\DeclareMathOperator{\Homka}{Hom_{\fld \mathrm{-alg}}}

\DeclareMathOperator{\Homabk}{Hom_{\Abk}}

\DeclareMathOperator{\Homfabk}{\underline{Hom}_{\Abk}}

\DeclareMathOperator{\K}{K}

\DeclareMathOperator{\Lie}{Lie}

\DeclareMathOperator{\R}{R}

\DeclareMathOperator{\Pic}{Pic}

\DeclareMathOperator{\Picf}{\underline{Pic}}

\DeclareMathOperator{\Specf}{\underline{Spec}}

\DeclareMathOperator{\Supp}{Supp}

\DeclareMathOperator{\Sym}{Sym}
\DeclareMathOperator{\Tg}{T}

\DeclareMathOperator{\Z}{Z}
\DeclareMathOperator{\Zero}{Z}

\DeclareMathOperator{\alb}{alb}

\DeclareMathOperator{\cl}{cl}
\DeclareMathOperator{\chr}{char}

\DeclareMathOperator{\dv}{div}
\DeclareMathOperator{\fml}{fml}

\DeclareMathOperator{\id}{id}
\DeclareMathOperator{\im}{im}

\DeclareMathOperator{\modu}{mod}

\DeclareMathOperator{\val}{v}

\newcommand{\Alba}[1]{\Alb(#1)}
\newcommand{\Fm}[2]{\fmlG_{#1,#2}}

\newcommand{\Um}[2]{\Upt\left(#1,#2\right)}
\newcommand{\Tm}[2]{\Trs\left(#1,#2\right)}
\newcommand{\Lm}[2]{L\left(#1,#2\right)}
\newcommand{\Jacc}[1]{\Jac\left(#1\right)}
\newcommand{\Jacm}[2]{\Jac\left(#1,#2\right)}
\newcommand{\Albm}[2]{\Alb(#1,#2)}

\newcommand{\AlbF}[1]{\Alb_{\fmlG}(#1)}

\newcommand{\AlbbF}[2]{\Alb_{\fmlG}^{(#1)}(#2)}

\newcommand{\Albfct}[1]{\Alb_{\mathrm{fct}}(#1)}
\newcommand{\Albmor}[1]{\Alb_{\mathrm{mor}}(#1)}

\newcommand{\Picor}[1]{\Pic^{0,\red}_{#1}}
\newcommand{\Picorf}[1]{\Picf^{0,\red}_{#1}}
\newcommand{\Decfc}[1]{\wh{\Decf_{#1}}}

\newcommand{\Divfc}[1]{\wh{\Divf_{#1}}}
\newcommand{\Divorc}[1]{\wh{\Divf^{0,\red}_{#1}}}
\newcommand{\Divor}[1]{\Divf^{0,\red}_{#1}}

\newcommand{\Divsyx}{\Divf_{Y/X}}

\newcommand{\Homsa}[1]{\mathrm{Hom}_{\,#1 \mathrm{-alg}}}
\newcommand{\Homs}[1]{\mathrm{Hom}_{#1}}

\renewcommand{\H}{\mathrm{H}}

\newcommand{\urq}[1]{\alb_{#1,#1_{\sing}}}

\newcommand{\urqqo}[1]{\alb_{\Xo,\Xo_{\sing}}^{(#1)}}
\newcommand{\Urq}[1]{\Alb(#1,#1_{\sing})}
\newcommand{\Urqq}[2]{\Alb^{(#1)}(#2,#2_{\sing})}

\newcommand{\WeilResC}[2]{\wh{\Pi}_{#1/#2}}

\newcommand{\ZO}[1]{\Z_0 (#1)}

\newcommand{\CHO}[1]{\CH_0 (#1)}
\newcommand{\CHOo}[1]{\CHO{#1}^{0}}
\newcommand{\CHSO}[1]{\CH_0 (#1,#1_{\sing})}
\newcommand{\CHSOo}[1]{\CHO{#1,#1_{\sing}}^{0}}
\newcommand{\CHm}[2]{\CH_0 (#1,#2)}
\newcommand{\CHmo}[2]{\CHm{#1}{#2}^{0}}

\newcommand{\Ab}{\mathsf{Ab}}
\newcommand{\Abk}{\mathcal{A} \mathit{b} / \fld}

\newcommand{\Abkp}{\mathsf{Ab} / \fld}

\newcommand{\AbVar}{\mathcal{A} \mathcal{V}}

\newcommand{\AlgS}[1]{\mathsf{Alg} / #1}
\newcommand{\Algk}{\mathsf{Alg} / \fld}

\newcommand{\Artk}{\mathsf{Art} / \fld}

\newcommand{\Gk}{\mathcal{G} / \fld}
\newcommand{\Gak}{\mathcal{G} \hspace{-.1em} \mathit{a} / \fld}
\newcommand{\aGk}{\mathit{a} \hspace{+.06em} \mathcal{G} / \fld}
\newcommand{\aGak}{\mathit{a} \hspace{+.06em} \mathcal{G} 
                              \hspace{-.1em} \mathit{a} / \fld}
\newcommand{\Gfk}{\mathcal{G} \hspace{-.01em} \mathit{f} / \fld}

\newcommand{\dGfk}{\mathit{d} \hspace{+.06em} \mathcal{G} 
                              \hspace{-.01em} \mathit{f} / \fld}

\newcommand{\Fctr}{\mathsf{Fctr}}

\newcommand{\Set}{\mathsf{Set}}

\newcommand{\Mr}{\mathsf{Mr}}

\newcommand{\Mrm}[2]{\Mr^{#1,#2}}

\newcommand{\MrCH}[1]{\Mr^{\CH}}
\newcommand{\MrCHm}[2]{\Mr^{\CHmo{#1}{#2}}}
\newcommand{\Morav}[1]{\mathsf{Mor} \pn{#1, \AbVar}}
\newcommand{\Fctav}[1]{\mathsf{Fct} \pn{#1, \AbVar}}

\newcommand{\Nat}{\mathbb{N}}
\newcommand{\Zint}{\mathbb{Z}}

\newcommand{\Prj}{\mathbb{P}}

\newcommand{\et}{\acute{\mathrm{e}} \mathrm{t}}

\newcommand{\mdl}{D}

\newcommand{\ord}{\mathrm{ord}}

\newcommand{\red}{\mathrm{red}}
\newcommand{\reg}{\mathrm{reg}}

\newcommand{\sing}{\mathrm{sing}}

\newcommand{\aug}{\eps}

\newcommand{\gal}{\sigma}

\newcommand{\hg}[1]{h_{#1}}

\newcommand{\homm}[1]{h^{(#1)}}
\newcommand{\emb}{\iota}

\newcommand{\morphism}{\psi}

\newcommand{\rmpC}{\varphi}

\newcommand{\rmp}[1]{\phe^{#1}}

\newcommand{\trafo}{\tau}

\newcommand{\Ampl}{\mathcal{A}}

\newcommand{\clfld}{\overline{\fld}}
\newcommand{\fld}{\mathit{k}}

\newcommand{\Algbr}{A}
\newcommand{\Algbrc}{\sA}

\newcommand{\Ring}{R}
\newcommand{\Rfin}{R}
\newcommand{\Sval}{S}

\newcommand{\Hpf}{\mathcal{H}}

\newcommand{\Pb}{\mathrm{P}}
\newcommand{\Gp}{\mathrm{G}}
\newcommand{\Fc}{F}
\newcommand{\Ft}{G}

\newcommand{\fmlG}{\mathcal{F}}

\newcommand{\Ga}{\mathbb{G}_{\mathrm{a}}}

\newcommand{\Gm}{\mathbb{G}_{\mathrm{m}}}

\newcommand{\Jac}{J}

\newcommand{\Lin}{\mathbb{L}}

\newcommand{\Trr}[1]{G^{(#1)}}
\newcommand{\Trs}{T}

\newcommand{\Trz}{\mathbb{T}}

\newcommand{\Upt}{U}

\newcommand{\Upf}{\mathbb{U}}

\newcommand{\Vsp}{V}

\newcommand{\tg}{t}
\newcommand{\pnt}{q}
\newcommand{\pntt}{p}

\newcommand{\Es}{E}

\newcommand{\Sing}{S}
\newcommand{\SgY}{\Sing_Y}

\newcommand{\Curv}{Z}
\newcommand{\Crvv}{C}
\newcommand{\Crv}{C}

\newcommand{\Xo}{X}
\newcommand{\Xb}{\overline{X}}

\newcommand{\pt}{\widetilde{p}}

\newcommand{\mc}{\widehat{\fm}}

\newcommand{\Oc}{\widehat{\sO}}

\newcommand{\Xt}{\widetilde{X}}

\newcommand{\Ct}{\widetilde{C}}

\newcommand{\Cy}{C^Y}
\newcommand{\Cyt}{\widetilde{\Cy}}

\newcommand{\Ld}{L^{\vee}}

\newcommand{\Il}{\mathcal{I}}

\newcommand{\alp}{\alpha}
\newcommand{\bet}{\beta}
\newcommand{\gam}{\gamma}
\newcommand{\del}{\delta}
\newcommand{\eps}{\varepsilon}

\newcommand{\lam}{\lambda}
\newcommand{\sig}{\sigma}

\newcommand{\phe}{\varphi}

\newcommand{\Gam}{\Gamma}

\newcommand{\Yps}{\Upsilon}

\newcommand{\fm}{\mathfrak{m}}

\newcommand{\fR}{\mathfrak{R}}

\newcommand{\sA}{\mathcal{A}}

\newcommand{\sD}{\mathcal{D}}

\newcommand{\sJ}{\mathcal{J}}
\newcommand{\sK}{\mathcal{K}}
\newcommand{\sL}{\mathcal{L}}

\newcommand{\sO}{\mathcal{O}}

\newcommand{\lrbt}[1]{\left[#1\right]}
\newcommand{\bigbt}[1]{\big[#1\big]}

\newcommand{\pn}[1]{(#1)}
\newcommand{\lrpn}[1]{\left(#1\right)}
\newcommand{\bigpn}[1]{\big(#1\big)}
\newcommand{\Bigpn}[1]{\Big(#1\Big)}

\newcommand{\bigst}[1]{\big\{#1\big\}}

\newcommand{\ra}{\rightarrow}

\newcommand{\dra}{\dashrightarrow}

\newcommand{\lra}{\longrightarrow}
\newcommand{\lla}{\longleftarrow}

\newcommand{\Llra}{\Longleftrightarrow}
\newcommand{\lmt}{\longmapsto}

\newcommand{\wt}[1]{\widetilde{#1}}
\newcommand{\wh}[1]{\widehat{#1}}
\newcommand{\ol}[1]{\overline{#1}}

\newcommand{\cut}{\cdot}
\newcommand{\isec}{\cap}
\newcommand{\dsum}{\bigoplus}

\newcommand{\tens}{\otimes}
\newcommand{\tensk}{\otimes_{\fld}}

\newcommand{\tensck}{\hspace{+.2em} \widehat{\otimes}_{\fld} \hspace{+.2em}}
\newcommand{\tms}{\times}

\newcommand{\bir}{\Yps}

\newcommand{\fis}{\#}
\newcommand{\blank}{\_}
\newcommand{\lul}{?}                    %%%  {\_}
\newcommand{\llul}{\hspace{+.05em} ?}   %%%  {\_}
\newcommand{\lull}{? \hspace{+.05em}}   %%%  {\_}
\newcommand{\see}{\textrm{see }}             %%%{\nearrow \,}
\newcommand{\seecite}{\textrm{see }}       %%%{\nearrow \,}
\newcommand{\laurin}{.}                     %%% {.} 
\newcommand{\laurink}{,}                   %%% {,} 
\newcommand{\vs}{12pt}
\newcommand{\vsvs}{24pt}

\newcommand{\bThm}{\begin{theorem}}
\newcommand{\eThm}{\end{theorem}}
\newcommand{\bAck}{\begin{acknowledgement}}
\newcommand{\eAck}{\end{acknowledgement}}
\newcommand{\bAlg}{\begin{algorithm}}
\newcommand{\eAlg}{\end{algorithm}}
\newcommand{\bAxm}{\begin{axiom}}
\newcommand{\eAxm}{\end{axiom}}
\newcommand{\bCas}{\begin{case}}
\newcommand{\eCas}{\end{case}}
\newcommand{\bClm}{\begin{claim}}
\newcommand{\eClm}{\end{claim}}
\newcommand{\bCcl}{\begin{conclusion}}
\newcommand{\eCcl}{\end{conclusion}}
\newcommand{\bCdn}{\begin{condition}}
\newcommand{\eCdn}{\end{condition}}
\newcommand{\bCjc}{\begin{conjecture}}
\newcommand{\eCjc}{\end{conjecture}}
\newcommand{\bCor}{\begin{corollary}}
\newcommand{\eCor}{\end{corollary}}
\newcommand{\bCrt}{\begin{criterion}}
\newcommand{\eCrt}{\end{criterion}}
\newcommand{\bDef}{\begin{definition}}
\newcommand{\eDef}{\end{definition}}
\newcommand{\bExm}{\begin{example}}
\newcommand{\eExm}{\end{example}}
\newcommand{\bExc}{\begin{exercise}}
\newcommand{\eExc}{\end{exercise}}
\newcommand{\bLem}{\begin{lemma}}
\newcommand{\eLem}{\end{lemma}}
\newcommand{\bNot}{\begin{notation}}
\newcommand{\eNot}{\end{notation}}
\newcommand{\bPrb}{\begin{problem}}
\newcommand{\ePrb}{\end{problem}}
\newcommand{\bPrp}{\begin{proposition}}
\newcommand{\ePrp}{\end{proposition}}
\newcommand{\bRmk}{\begin{remark}}
\newcommand{\eRmk}{\end{remark}}
\newcommand{\bSol}{\begin{solution}}
\newcommand{\eSol}{\end{solution}}
\newcommand{\bSmr}{\begin{summary}}
\newcommand{\eSmr}{\end{summary}}
\newcommand{\bVar}{\begin{variant}}
\newcommand{\eVar}{\end{variant}}
\newcommand{\bPf }{\begin{prooof}}
\newcommand{\ePf }{\end{prooof}}
\theoremstyle{plain}
\newtheorem{theorem}{Theorem}[section]
\newtheorem{axiom}[theorem]{Axiom}
\newtheorem{conjecture}[theorem]{Conjecture}
\newtheorem{corollary}[theorem]{Corollary}
\newtheorem{criterion}[theorem]{Criterion}
\newtheorem{lemma}[theorem]{Lemma}
\newtheorem{problem}[theorem]{Problem}
\newtheorem{proposition}[theorem]{Proposition}

\theoremstyle{definition}
\newtheorem{acknowledgement}[theorem]{Acknowledgement}
\newtheorem{algorithm}[theorem]{Algorithm}
\newtheorem{case}[theorem]{Case}
\newtheorem{claim}[theorem]{Claim}
\newtheorem{condition}[theorem]{Condition}
\newtheorem{conclusion}[theorem]{Conclusion}
\newtheorem{definition}[theorem]{Definition}
\newtheorem{example}[theorem]{Example}
\newtheorem{exercise}[theorem]{Exercise}
\newtheorem{notation}[theorem]{Notation}
\newtheorem{remark}[theorem]{Remark}
\newtheorem{solution}[theorem]{Solution}
\newtheorem{summary}[theorem]{Summary}

\newenvironment{prooof}[1][Proof]{\textbf{#1.} }
{\ \rule{0.5em}{0.5em}} 
\newenvironment{variant}[1][Variant]{\textbf{#1.} }
{\ \rule{0.5em}{0.5em}}

\begin{document}

\centerline{ } 
\vspace{35pt} 
\centerline{\LARGE{Albanese Varieties of Singular Varieties}} 
\vspace{4pt} 
\centerline{\LARGE{over a Perfect Field}} 
\vspace{22pt}
\centerline{\large{Henrik Russell}} 
%\vspace{10pt}
%%\centerline{\large{November 2010}} 
\vspace{22pt}

%\title{Albanese of Singular Varieties\\ 
%         over a Perfect Field} 
%\author{Henrik Russell}
%\maketitle

%%%%%%%%  Abstract  %%%%%%%%%%%%%%%%%%%%%%%%%%%%%%%%%%%%%%%%%%%%%%%%%%%%%%%%%
\begin{abstract}
Let $X$ be a projective variety, possibly singular. 
A generalized Albanese variety of $X$ was constructed by 
Esnault, Srinivas and Viehweg 
for algebraically closed base field as a universal regular quotient 
of a relative Chow group of 0-cycles of degree 0 modulo rational equivalence. 
In this paper,  
we obtain a functorial description of the Albanese of Esnault-Srinivas-Viehweg 
over a perfect base field, 
using duality theory of 1-motives with unipotent part. 
\end{abstract}

%\newpage 
%\quad 
%\newpage 

\tableofcontents{} 
%\newpage 
%\quad 
%\newpage 

\setcounter{section}{-1}

%%%%%%%%%%%%  Introduction  %%%%%%%%%%%%%%%%%%%%%%%%%%%%%%%
\section{Introduction}

With a smooth projective variety $X$ over a field $k$ 
there is an associated Albanese variety $\Alba{X}$, 
defined by a universal mapping property 
for rational maps from $X$ to abelian varieties. 
Equivalently, $\Alba{X}$ can be considered as a universal quotient 
of the Chow group $\CHOo{X}$ 
of 0-cycles of degree 0 on $X$ modulo rational equivalence. 
This second characterization was used by 
Esnault, Srinivas and Viehweg to define an Albanese variety $\Urq{X}$ 
for possibly singular $X$ over an algebraically closed field $k$. 
Here $\CHO{X}$ was replaced by the relative Chow group of 0-cycles 
$\CHSO{X}$ in the sense of Levine-Weibel \cite{LW}, 
which coincides with $\CHO{X}$ if $X_{\sing} = \varnothing$. 
Then $\Urq{X}$ is characterized by a universal mapping property 
for rational maps from $X$ to smooth connected commutative algebraic groups 
that factor through $\CHSOo{X}$. 

The aim of this paper is to give a functorial description of the Albanese $\Urq{X}$ 
of Esnault, Srinivas and Viehweg for possibly singular $X$ over a perfect field 
of arbitrary characteristic, 
see Theorems \ref{main_result}, \ref{main_result_2} and \ref{descentUrq}. 
For this purpose we use the concept of generalized Albanese varieties 
over a perfect field and duality of 1-motives with unipotent part from \cite{Ru2}. 
In this way we obtain a functor dual to $\Urq{X}$ as well. 
The construction is explicit. 

The case that the base field $k$ is of characteristic 0 and algebraically closed 
was done in \cite{Ru}. 
While the way of procedure appears analogous in parts, 
it was necessary to involve different methods, 
e.g.\ formal groups in positive characteristic behave distinctly 
from the characteristic 0 case 
and we cannot assume in general a resolution of singularities.

%\newpage 

\subsection{Leitfaden} 

We make use of the notions described in \cite{Ru2}, in particular 
\emph{categories of rational maps from a variety $X$ to algebraic groups} 
\cite[Definition~2.7]{Ru2} %Definition~\ref{CatMr} 
and \emph{universal objects} for such categories of rational maps 
\cite[Definition~2.11]{Ru2}. %Definition~\ref{univObj}. 
%\medskip 

Let $X$ be a (singular) projective variety over an algebraically closed field $k$. 
The Albanese variety of Esnault-Srinivas-Viehweg  $\Urq{X}$ 
is by definition the universal object for the category $\MrCH{X}$ 
of rational maps from $X$ to smooth connected commutative algebraic groups 
factoring through $\CHSOo{X}$ 
($\see$Definition~\ref{MrCH(X)}), 
therefore it is also called \emph{the universal regular quotient of $\CHSOo{X}$}. 
Suppose $\pi: Y \ra X$ is a projective resolution of singularities 
(we will later replace this by a morphism 
that is always available in any characteristic). 
%Via pull-back we may consider a rational map $\phe:X \dra G$ 
%as a rational map from $Y$ to $G$. 
We give a characterization of rational maps factoring through $\CHSOo{X}$ 
in terms of a formal subgroup of the functor $\Divf_Y$ 
of relative Cartier divisors on $Y$: 

A rational map $\phe: X \dra G$ to a smooth connected commutative 
algebraic group induces a transformation $\trafo_{\phe}: \Ld \lra \Divf_Y$ 
($\seecite$\cite[Definition~2.4]{Ru2}), %($see$Definition~\ref~{induced Trafo}) 
where $L$ is the affine part of $G$ and $\Ld$ the Cartier dual. 
Let $\Mr_{\Divor{Y/X}}$ be the category of all those rational maps $\phe$ 
from $Y$ to smooth connected commutative algebraic groups 
such that the image of the induced transformation $\trafo_{\phe}$ 
is contained in $\Divor{Y/X}$, 
``the kernel of $\Divf_Y$ under push-forward from $Y$ to $X$''. 
The precise definition of $\Divor{Y/X}$ 
(see Definitions \ref{Div_Z/C^0} and \ref{Div_Y/X^0}) 
%and \ref{Div_Y/X^0,S}) 
requires some technical issues: 
so-called \emph{formal divisors} (see Definition~\ref{FDiv}) 
and restriction to Cartier curves. 
Then the category $\MrCH{X}$ is equal to $\Mr_{\Divor{Y/X}}$.

Using the result \cite[Theorem~2.12]{Ru2} 
% Theorem~\ref{Exist univObj} 
about the existence of universal objects 
for categories of rational maps, %to algebraic groups, 
we can re-prove the existence 
of the universal regular quotient $\Urq{X}$ of $\CHSOo{X}$ 
(alternative to \cite{ESV}) 
and obtain an explicit and functorial construction 
(cf.\ Section~\ref{subsubsec:Exist_UnivRegQuot}): 

\bThm 
\label{main_result}
Let $X$ be a projective variety over an algebraically closed field $k$, 
and let $\Xt \ra X$ be a projective resolution of singularities. 
Then the universal regular quotient $\Urq{X}$ of $\CHSOo{X}$ exists 
(as an algebraic group) 
and its dual (in the sense of 1-motives) represents the functor 
\[ \Divor{\Xt/X} \lra \Picorf{\Xt} 
\] 
that is the natural transformation of functors 
which assigns to a relative divisor the class of its associated line bundle. 
$\Picor{\Xt}$ is represented by an abelian variety and
$\Divor{\Xt/X}$ by a dual-algebraic formal group 
(i.e.\ a formal group whose Cartier dual is an algebraic affine group). 
\eThm 

In the general case, i.e.\ if we do not assume a resolution of singularities, 
we can perform the same procedure, replacing \emph{resolution of singularities} 
by a suitable blowing up of $X$. 
The construction is elementary, 
in particular it does not use \emph{alterations}. 
We obtain (cf.\ Section~\ref{no Desingularization}): 

\bThm 
\label{main_result_2}
Let $X$ be a projective variety over an algebraically closed field $k$. 
Suppose $\pi: Y \lra X$ is a birational projective morphism 
with the property that every rational map from $X$ to an abelian variety $A$ 
extends to a morphism from $Y$ to $A$, 
and $\Picor{Y}$ is an abelian variety. 
Then the functor $\bigbt{\Divor{Y/X} \lra \Picorf{Y}}$ 
is dual (in the sense of 1-motives) to 
the universal regular quotient $\Urq{X}$ of $\CHSOo{X}$. 
A morphism $\pi: Y \lra X$ with the required properties does always exist 
and can be constructed as a blowing up of an ideal sheaf on $X$. 
In particular, $\Urq{X}$ exists (as an algebraic group). 
\eThm 

A descent of the base field yields 
(cf.\ Theorem~\ref{descent_UnivRegQuot}): 

\bThm 
\label{descentUrq}
Let $X$ be a projective irreducible variety defined over a perfect field $\fld$. 
Let $\clfld$ be an algebraic closure of $\fld$. 
Then $\Alb\bigpn{X \tens_{\fld} \clfld,\pn{X \tens_{\fld} \clfld}_{\sing}}$ 
is defined over $\fld$. 
\eThm 

Some hint for the reader: as the notations from \cite{Ru2} are used 
throughout this paper, it might be helpful to know that there is a 
Glossary of Notations on page \pageref{Glossar}. 

I am very thankful to H\'el\`ene Esnault, Kazuya Kato and Andre Chatzistamatiou 
for discussions and advices.

\newpage 
%\quad 
%\newpage 

\section{Prerequisites} 
\label{prerequisites}

Throughout this note the following assumptions hold 
(unless stated otherwise). 
$k$ is an algebraically closed field (of arbitrary characteristic), 
except in No.~\ref{UnivRegQuot_descent}, 
where $k$ is only assumed to be perfect. 
$X$ is a projective variety over $k$ 
and $\pi:Y\lra X$ a resolution of singularities, 
i.e.\ $Y$ is regular and $\pi$ is a projective birational morphism 
which is an isomorphism over the regular points of $X$. 
In No.s~\ref{no Desingularization} and \ref{UnivRegQuot_descent}
we do not assume the existence of a resolution of singularities. 
$U\subset X_{\reg}$ is an open dense subset of $X$, 
contained in the regular locus of $X$. 
%As $\pi$ is an isomorphism on $\pi^{-1}U$, 
%we identify $U$ with its inverse image in $Y$. 
We consider the category $\MrCH{X}$ of morphisms 
$\phe:U\lra G$ from $U$ to algebraic groups $G$ 
factoring through $\CHSOo{X}$ ($\see$Definition~\ref{MrCH(X)}). 
Here we assume algebraic groups $G$ 
always to be smooth, connected and commutative, 
unless stated otherwise. 

\subsection{Chow Group of Points on Singular Varieties} 
\label{Chow_LW}

We recall the relative Chow group $\CHSOo{X}$
of 0-cycles of degree 0 modulo rational equivalence from \cite{LW}. 
%see also \cite{ESV}, \cite{BiS}. 

\bDef 
\label{Cartier-curve} 
Let $S$ be a closed proper subset of $X$. \\
A \emph{Cartier curve} in $X$ relative to $S$ 
is a curve $C\subset X$ satisfying 

\begin{tabular}{rl}
(a) & $C$ is pure of dimension 1\laurin \\
(b) & No component of $C$ is contained in $S$\laurin \\
(c) & If $p\in S$ is a point of $C\cap S$, 
         then the ideal of $C$ in $\sO_{X,p}$ 
         is gene- \\ 
     & rated by a regular sequence 
         (consisting of $\dim X - 1$ elements)\laurin \\ 
\end{tabular} 
\eDef 

\bDef 
\label{R(C,X)} 
Let $C$ be a Cartier curve in $X$ relative to $X\setminus U$. 
Let $\gamma_1,\ldots, \gamma_r$ be the generic points of $C$. 
%$Z\in\Cp\pn{C}$, 
%where $\Cp\pn{C}$ is the set of the irreducible components of $C$. 
Let $\sO_{C,\Theta}$ be the semilocal ring on $C$ at the points in 
$\Theta=\lrpn{C\setminus U} \cup \lrpn{ \gamma_1,\ldots,\gamma_r} $. 
%There is a natural map on unit groups 
%\[ \vartheta_{C,U}:\sO_{C,\Theta}^{*}\lra\bigoplus_{Z\in\Cp\pn{C}}\K_{Z}^{*} 
%\] 
Then define 
\[ \K\lrpn{C,U}^{*} =\, \sO_{C,\Theta}^{*} \laurin %\im\vartheta_{C,U} 
\]
\eDef 

\bDef 
\label{div(f)_C} 
Let $C$ be a Cartier curve in $X$ relative to $X\setminus U$ 
and let $\nu:\Ct\lra C$ be its normalization. 
For $f\in\K\lrpn{C,U}^{*}$ and $p\in C(k)$ let 
\[ \ord_{p}\lrpn{f}=\sum_{\pt\ra p}\val_{\pt}\lrpn{\widetilde{f}} 
\] 
where $\widetilde{f} := \nu^{\fis}f$ %$ = f \circ \nu \in \sK_{\Ct}$ 
and $\val_{\pt}$ is the discrete valuation attached to the point $\pt\in\Ct(k)$ 
above $p\in C(k)$ 
(cf.\ \cite[Exm.~A.3.1]{F}). 
Define the divisor of $f$ to be 
\[ \dv\pn{f}_{C}=\sum_{p\in C}\ord_{p}\lrpn{f}\,\left[p\right] \laurin 
\] 
\eDef  

Immediately from Definition~\ref{div(f)_C} we get 

\bLem 
\label{normlzd} 
Let $C$ be a Cartier curve in $X$ relative to $X\setminus U$ 
with normalization $\nu:\Ct\lra C$. 
If $\rmpC:C\cap U\lra G$ is a morphism
from $C\cap U$ to an algebraic group $G$, then for each 
$f\in\K\lrpn{C,U}^*$ it holds 
(with \,$\phe\bigpn{\sum n_i \,p_i} := \sum n_i \, \phe\pn{p_i}$) 
\[ \rmpC\lrpn{\dv\pn{f}_{C}} = 
   \lrpn{\rmpC\circ\nu}\lrpn{\dv\pn{f\circ\nu}_{\Ct}} %\laurin 
\] 
%where \;$\phe\bigpn{\sum n_i \,p_i} := \sum n_i \, \phe\pn{p_i}$. 
\eLem 

%\bPf  
%By Definition~\ref{div(f)_C} we have 
%\begin{eqnarray*}
%\rmpC\lrpn{\dv\pn{f}_{C}} 
% & = & \sum_{p\in C}\ord_{p}\lrpn{f}\,\rmpC(p)\\
% & = & \sum_{p\in C}\sum_{\pt\ra p}\val_{\pt}\lrpn{f\circ\nu}\,\rmpC(p) \\
% & = & \sum_{\pt\in\Ct}\val_{\pt}\lrpn{f\circ\nu}\, 
%       \lrpn{\rmpC\circ\nu}\lrpn{\pt}\\
% & = & \lrpn{\rmpC\circ\nu}\lrpn{\dv\pn{f\circ\nu}_{\Ct}}
%\end{eqnarray*} 
%\hfill 
%\ePf 

\bDef 
\label{CH_0(X)_deg0}
Let $\ZO{U}$ be the group
of 0-cycles on $U$, set 
\[ \fR_0\pn{X,U}=\left\{ \lrpn{C,f}\left|\begin{array}{c}
C\textrm{ is a Cartier curve in }X\textrm{ relative to }X\setminus U\\
\textrm{and }f\in\K\lrpn{C,U}^{*}\end{array}\right.\right\} 
\] 
and let $\R_{0}\pn{X,U}$ be the subgroup of $\ZO{U}$
generated by elements $\dv\pn{f}_{C}$ with 
$\lrpn{C,f}\in\fR_0\pn{X,U}$.
Then define 
\[ \CHSO{X}=\ZO{U}/\R_{0}\pn{X,U} \laurin 
\] 
Let $\CHSOo{X}$ be the subgroup of $\CHSO{X}$
of cycles $\zeta$ with \\ 
$\deg\zeta|_{W}=0$ for all irreducible components $W$ of $U$. 
%$W\in\Cp\lrpn{U}$
\eDef  

\bRmk  
The definition of $\CHSO{X}$ and $\CHSOo{X}$
is independent of the choice of the dense open subscheme $U\subset X_{\reg}$
($\seecite$\cite[Corollary 1.4]{ESV}). 
\eRmk  

\bRmk  
By our terminology a curve is always reduced, 
in particular a Cartier curve. 
Allowing non-reduced Cartier curves 
(which is the common definition) 
does not change the groups $\CHSO{X}$ and $\CHSOo{X}$, 
see \cite[Lemma 1.3]{ESV} for more explanation. 
\eRmk 

\bRmk 
The definition of $\CHSO{X}$ is made in such a way 
that for any (possibly singular) curve $\Crv$ 
the Chow group of points coincides with the Picard group: 
$\CHSO{\Crv} = \Pic\pn{\Crv}$ and 
$\CHSOo{\Crv} = \Pic^0\pn{\Crv}$, 
see \cite[Proposition~1.4]{LW}. 
\eRmk

\subsection{Functors of Relative Cartier Divisors} 
\label{FuncDiv}

Relative Cartier divisors will be used 
for the description of (the duals of) generalized Albanese varieties. 
For the sake of functoriality we need not only to consider Cartier divisors on $X$, 
but also their infinitesimal deformations. 
Of great relevance will be the completion $\Divfc{X}$ 
of the functor $\Divf_X$ of \emph{relative Cartier divisors on $X$} 
(see \cite[No.~2.1]{Ru}), 
which assigns to a finite (i.e.\ finite dimensional) $k$-algebra $R$ 
the abelian group 
\[ \Divfc{X}(R) = 
   \Gam\bigpn{X \tens_k R, \pn{\sK_X \tens_k R}^* / \pn{\sO_X \tens_k R}^*} 
\]
where the star ${}^*$ denotes the unit groups. 
The functor \;$\Divfc{X}: \Artk \lra \Ab$\; 
from the category of finite $k$-algebras to the category of abelian groups 
is a formal group ($\seecite$\cite[Proposition 2.1]{Ru2}). 
%($\see$\ref{Divf formal group}). 

Let $\morphism: V \lra X$ be a morphism of varieties. 

\bPrp 
\label{Decf formal group} 
Let $\Decfc{X,V}$ (cf.\ \cite[Definition 2.11]{Ru}) 
be the subfunctor of $\Divfc{X}$ 
consisting of those relative Cartier divisors on $X$ 
whose support ($\seecite$\cite[Definition 2.2]{Ru2}) 
%($\see$Definition~\ref{Support}) 
intersects $\morphism\pn{V}$ properly. 
Then 
\;$\Decfc{X,V} : \Artk \lra \Ab$\; 
is a formal $k$-group. 
\ePrp 

\bPf 
Analogous to \cite[Proposition 2.1]{Ru2}. 
\ePf 

\bPrp 
\label{_.V} 
The pull-back of Cartier divisors from $X$ to $V$ 
induces a homomorphism of formal groups 
\quad$ \lul\cut V: \Decfc{X,V} \lra \Divfc{V} \laurin $ 
%of relative Cartier divisors from $X$ to $V$ 
\ePrp 

\bPf 
Straightforward. 
\ePf 

%(Cf.\ \cite[Definition 2.12]{Ru}.) 

%\bPrp 
%\label{Div_X^S formal group} 
%If $S$ is a closed subscheme of $X$, 
%let $\Divsfc{X}{S}$ be the subfunctor of $\Divfc{X}$ 
%consisting of those relative Cartier divisors on $X$ 
%whose support is contained in $S$. 
%Then \;$\Divsfc{X}{S} : \Artk \lra \Ab$\; 
%is a formal $k$-group. 
%\ePrp 
%
%\bPf 
%Analogous to \cite[Proposition 2.1]{Ru2}. 
%\ePf 

%\newpage 

\section{Formal Divisors} 
\label{FormalDivisors}

For Cartier divisors on a $k$-scheme $X$ 
there is a notion of pull-back, 
but not a push-forward. %of relative Cartier divisors. 
If $X$ is a normal scheme, 
the group of Cartier divisors on $X$ 
can be identified with the subgroup of locally principal Weil divisors, 
and there is a push-forward of Weil divisors. 
But for Weil divisors it does not make sense to speak about 
infinitesimal deformations, because Weil divisors are formal sums 
of prime divisors and prime divisors are always reduced. 

We introduce a formal group functor $\FDivf_X$ of formal divisors on $X$ 
that admits a push-forward and a natural transformation 
\;$\fml: \Divfc{X} \lra \FDivf_X$. 
The group of formal infinitesimal divisors $\LDiv\pn{X}$ 
from \cite[No.~3.3]{Ru} 
is the Lie algebra of $\FDivf_X$. 

In this section we consider the case that $X$ is a curve $\Curv$.

%\newpage 

\subsection*{Functor of Formal Divisors} 

Let $\Curv$ be a curve over $k$. 
%over $k$ \\  (an algebraically closed field of characteristic 0). 

\bDef 
\label{FDiv} 
Define the formal $k$-group functor of \emph{formal divisors on $\Curv$} 
\[ \FDivf_{\Curv} : \Artk \lra \Ab 
\] 
from the category of finite $k$-algebras 
to the category of abelian groups  by 
\[ \FDivf_{\Curv} = \dsum_{\pnt\in\Curv(k)} 
   \Homfabk\lrpn{\Oc_{\Curv,\pnt}^*,k^*} 
\] 
where $\Homfabk$ denotes the $k$-group functor 
of homomorphisms of $k$-group functors, 
%($\see$begin of Section~\ref{Cartier-Duality}), 
and $\sA^*$ stands for the $k$-group functor 
$\Gm\lrpn{\llul \tensck \sA}$ for any profinite $k$-algebra $\sA$. \\ 
As Weil restriction is right-adjoint to base extension, 
cf.\ \cite[1.1.2]{Ru2}, %Proposition~\ref{adjoint functors}, 
the $\Ring$-valued points of $\FDivf_{\Curv}$ are given by 
\[ \FDivf_{\Curv}(\Ring) = \dsum_{\pnt\in\Curv(k)} 
   \Homabk\Big(\Gm\big(\llul\tensck\Oc_{\Curv,\pnt}\big),
               \Gm\big(\llul\tens_k\Ring\big)\Big) \laurin 
\] 
\eDef 

\bPrp 
\label{FDiv repr}
The functor $\FDivf_{\Curv}$ of formal divisors on $\Curv$ 
is a formal $k$-group. 
\ePrp 

\bPf 
The formal $k$-group functor 
\begin{eqnarray*}
\Homfabk\bigpn{\Oc_{\Curv,\pnt}^*,k^*} 
% & = & \Homabk\bigpn{\Oc_{\Curv,\pnt}^*,\Lin_{\lul}} \\ 
 & = & \Homabk\bigpn{\Oc_{\Curv,\pnt}^*,\lull} \circ \Gm \circ \pn{\blank \tens_k \lull}  
\end{eqnarray*}
as a composition of left-exact functors (the tensor product is over a field), 
is left-exact, hence a formal group ($\seecite$\cite[Prop.~1.7]{Ru2}). 
%($\see$Proposition \ref{formalScheme}). 
Finite direct sums of formal groups are formal groups. 
The finite subsets $\Es \subset \Curv(k)$ form a directed inductive system, 
ordered by inclusion. 
Thus 
\[ \FDivf_{\Curv} = \varinjlim_{\Es\subset\Curv(k)} \; \dsum_{\pnt\in\Es} \;
                                \Homfabk\lrpn{\Oc_{\Curv,\pnt}^*,k^*} 
\] 
is the limit of a directed inductive system of formal groups, 
hence a formal group. 
\ePf 

\bRmk 
\label{Oc affine}
The $k$-group functor 
$\Oc_{\Curv,\pnt}^* = \Gm\bigpn{\llul\tensck\Oc_{\Curv,\pnt}}$ 
is an affine group 
by Lemma~\ref{Lin_profiniteRing} below, 
since the completion $\Oc_{\Curv,\pnt}$ of the local ring at $\pnt\in\Curv$ 
is a profinite $k$-algebra. 
%This is neither used nor implied by Proposition \ref{FDiv repr}. 
\eRmk 

\bLem 
\label{Lin_profiniteRing} 
Let $\fld$ be a ring and 
$\Algbrc$ a profinite $\fld$-algebra 
which is a projective $\fld$-module. 
Let $\Ft$ be an $\Algbrc$-functor 
which is represented by an affine $\Algbrc$-scheme. 
Then the induced $\fld$-functor \;$\Ft\lrpn{\llul\tensck\Algbrc}$\; 
is represented by an affine $\fld$-scheme. 
\eLem 

\bPf 
We adapt the proof of \cite[I, \S~1, Proposition~6.6]{DG} 
to the situation where $\Algbrc$ is, 
instead of a finite $\fld$-algebra, 
a profinite $\fld$-algebra which is a projective $\fld$-module. 
%(cf.\ Remark~\ref{WeilRes_profinite}). 

First assume that $\Ft = \Specf_{\Algbrc} \Sval$, where 
$\Sval = \Sym_{\Algbrc} \lrpn{\Vsp\tensk\Algbrc}$ 
is the symmetric algebra of a $\fld$-module of the form 
$\Vsp \tensk \Algbrc$. 
Let $\{\Algbr_i\}$ be a projective system of discrete finite quotients of 
$\Algbrc$ such that $\Algbrc = \varprojlim \Algbr_i$, 
and let $\Hpf = \Homck(\Algbrc,\fld)$. 

For $\Ring \in \Algk$ 
\begin{eqnarray*}
\Ft \lrpn{ \Ring \tensck \Algbrc } 
 & = & \Homsa{\Algbrc} \lrpn{ \Sval, \Ring \tensck \Algbrc } \\
 & = & \Homs{\Algbrc} \lrpn{ \Vsp \tensk \Algbrc, \Ring \tensck \Algbrc } \\
% & = & \Homk \lrpn{ \Vsp, \Ring \tensck \Algbrc } \\
 & = & \Homk \lrpn{ \Vsp, \Ring \tensck \varprojlim \Algbr_i } \\
 & = & \Homk \lrpn{ \Vsp, \varprojlim \lrpn{\Ring \tensk \Algbr_i } } \\
 & = & \varprojlim \; \Homk \lrpn{ \Vsp, \Ring \tensk \Algbr_i } \\
 & = & \varprojlim \; \Homk \lrpn{ \Vsp \tensk \Homk(\Algbr_i,\fld), \Ring } \\
 & = & \Homk \lrpn{ \varinjlim \lrpn{\Vsp \tensk \Homk(\Algbr_i,\fld)}, \Ring } \\
 & = & \Homk \lrpn{ \Vsp \tensk \varinjlim \; \Homk(\Algbr_i,\fld), \Ring } \\
 & = & \Homk \lrpn{ \Vsp \tensk \Homck(\varprojlim \Algbr_i,\fld), \Ring } \\
 & = & \Homk \lrpn{ \Vsp \tensk \Homck(\Algbrc,\fld), \Ring } \\
 & = & \Homka \big( \Sym_{\fld}(\Vsp \tensk \Hpf), \Ring \big) 
\end{eqnarray*}
where we used in the fourth equation 
\begin{eqnarray*}
\Ring \tensck \varprojlim \Algbr_i & = & \varprojlim \big(\Ring \tensk \Algbr_i \big) \\
 & = & \left\{ \left. \sum_i r_i \tens a_i \; \right|  
            \begin{array}{l}
            r_i \in \Ring, \quad \; a_i \in \Algbrc \;\; \textrm{ s.t.} \\
            \ol{a_i} \in \ker \big( A_{i} \ra A_{j} \big) \;\;
            \textrm{for} \;\; i > j 
            \end{array}
            \right\}
\end{eqnarray*} 
and in the eighth equation, according to 
\cite[Chapter 2, Exercise 20, p.~33]{A} 
\[  \Vsp \tensk \varinjlim H_i \; = \; \varinjlim \big(\Vsp \tensk H_i\big) \laurin 
\] 
Thus one has 
\[  \WeilResC{\Algbrc}{\fld} \;
    \Specf_{\Algbrc} \big( \Sym_{\Algbrc} (\Vsp\tensk\Algbrc) \big) = 
    \Specf_{\fld} \big( \Sym_{\fld} (\Vsp\tensk\Hpf) \big)
\] 
where 
\;$\WeilResC{\Algbrc}{\fld}: \Fctr(\AlgS{\Algbrc},\Set) \lra \Fctr(\Algk,\Set)$\; 
is the functor given by \;$\Ft \lmt \Ft \lrpn{ \llul \tensck \Algbrc }$. 

For the case $\Ft = \Specf_{\Algbrc} \Sval$, 
where $\Sval$ is an arbitrary $\Algbrc$-algebra, 
let $\Il$ be the kernel of the canonical homomorphism 
$\Sym_{\Algbrc} (\Sval \tensk \Algbrc) \lra \Sval$. 
Then $\Sval$ is identified in $\AlgS{\Algbrc}$ 
with the amalgamated sum due to the diagram 
\[  \Sym_{\Algbrc} (\Sval \tensk \Algbrc)  \overset{\emb}\lla 
    \Sym_{\Algbrc} (\Il \tensk \Algbrc)  \overset{\aug}\lra \Algbrc 
\] 
%\[  \xymatrix{ 
%    \Sym_{\Algbrc} (\Sval \tensk \Algbrc)  & &  \Algbrc  \\
%    & \Sym_{\Algbrc} (\Il \tensk \Algbrc) \ar[ul]^{\emb} \ar[ur]_{\aug} & 
%                   }
%\] 
where $\emb$ is the canonical homomorphism and 
$\aug(\Il \tensk \Algbrc) = 0$. 
The functor $\WeilResC{\Algbrc}{\fld}$ is left-exact, 
i.e.\ it commutes with finite projective limits. 
%due to the definition of projective limits in $\Fctr(\AlgS{\Algbrc},\Set)$. 
Then $\WeilResC{\Algbrc}{\fld} \Specf_{\Algbrc} (\Sval)$ is the fibre product 
obtained from the transformation of the amalgamated sum by 
$\WeilResC{\Algbrc}{\fld} \circ \Specf_{\Algbrc}$. 
Thus $\WeilResC{\Algbrc}{\fld} \Ft$ is represented 
by an affine $\fld$-scheme, 
since it is given by the fibre product of affine $\fld$-schemes. 
\ePf 

\bDef 
\label{phi_*}
Let $\pi:\Curv\lra\Crvv$ be a finite morphism 
%of degree 1 ($\see$Definition~\ref{deg1})
of curves over $k$. 
The \emph{push-forward} of formal divisors 
\[ \pi_* : \FDivf_{\Curv} \lra \FDivf_{\Crvv} 
\] 
is the homomorphism of formal $k$-group functors 
induced by the homomorphisms 
\begin{eqnarray*}
\Homfabk\lrpn{\Oc_{\Curv,\pnt}^*,k^*} & \lra & 
\Homfabk\lrpn{\Oc_{\Crvv,\pi(\pnt)}^*,k^*} \\
h & \lmt & h \circ \Bigpn{\wh{\pi^{\fis}_{\pnt}}}^* 
\end{eqnarray*} 
where $\pnt\in\Curv(k)$ 
and $\wh{\pi^{\fis}_{\pnt}}: \Oc_{\Crvv,\pi(\pnt)} \lra \Oc_{\Curv,\pnt}$ 
is the local homomorphism of complete local rings associated with $\pi$, 
giving 
$\Bigpn{\wh{\pi^{\fis}_{\pnt}}}^* = \Gm\Bigpn{\llul\tensck\wh{\pi^{\fis}_{\pnt}}}$. 
\eDef 

%\bPf Evident. \ePf 

\bPrp 
\label{fml} 
Let $\Curv$ be a normal curve over $k$. 
There is a homomorphism of formal $k$-groups 
%($\see$Remark~\ref{fmlGroup_on_finiteAlgebras})
\[ \fml : \Divfc{\Curv} \lra \FDivf_{\Curv} 
\] 
for each finite $k$-algebra $R$ given by 
\begin{eqnarray*} 
\fml\pn{R} \; : \; 
\Gamma\Big(\Curv\tens R,\left.\sK_{\Curv\tens R/R}^*\right/\sO_{\Curv\tens R}^*\Big) 
 & \lra & \dsum_{\pnt\in\Curv(k)} \Homabk\lrpn{\Oc_{\Curv,\pnt}^*,R^*} \\ 
\sD & \lmt & \sum_{\pnt\in\Curv(k)} (\sD,\lull)_{\pnt} 
\end{eqnarray*} 
where $\lrpn{\rmpC,\lull}_{\lul}:\sK_{\Curv}^{*}\times\Curv\lra G$ 
is the local symbol associated with a rational map $\rmpC$ 
from $\Curv$ to an algebraic group $G$. 
\ePrp 

\bPf 
The first task is to show that for $\sD\in\Divfc{\Curv}(R)$ 
and $\pnt\in\Curv(k)$ the expression 
$(\sD,\lull)_{\pnt} : \Oc_{\Curv,\pnt}^* \lra R^*$ 
yields a well defined homomorphism of abelian groups. 
After that the naturality will be shown. 

Let $\del\in\Gamma\big(U_{\pnt}\tens R,\sK_{\Curv\tens R/R}^*\big)$ 
be a local section representing $\sD$ in a neighbourhood 
$U_{\pnt}\subset\Curv$ of $\pnt$. The restriction of $\del$ to a 
suitable open dense subset $U\subset U_{\pnt}$ gives an element of 
$\Gamma\big(U\tens R,\sO_{\Curv\tens R}^*\big) 
 = \Gm\bigpn{\sO_{\Curv}(U)\tens_k R} = \Lin_R\bigpn{\sO_{\Curv}(U)}$. 
Thus $\del$ can be seen as a rational map $\del: \Curv \dra \Lin_R$. 

Let $\del'\in\Gamma\big(U_{\pnt}'\tens R,\sK_{\Curv\tens R/R}^*\big)$ 
be another local section representing $\sD$ in a neighbourhood 
$U_{\pnt}'\subset\Curv$ of $\pnt$. 
Then $\del$ and $\del'$ differ at $\pnt$ by a unit 
$\gam\in\Gamma\lrpn{U_{\pnt}''\tens R,\sO_{\Curv\tens R}^*}$, 
i.e.\ $\del'=\del\cdot\gam$. 
The corresponding rational map of this unit $\gam: \Curv \dra \Lin_R$ 
is regular at $\pnt$. 
Therefore for each $f\in\Oc_{\Curv,q}^*$ we compute 
$(\gam,f)_{\pnt} = \gam(\pnt)^{\val_{\pnt}(f)} = 1$ 
since $\val_{\pnt}(f) = 0$. 
By \cite[Proposition 3.15]{Ru} %Proposition~\ref{locSymb_additiv} 
one obtains 
$(\del',f)_{\pnt} = (\del\gam,f)_{\pnt} 
                = (\del,f)_{\pnt}\cdot(\gam,f)_{\pnt} 
                = (\del,f)_{\pnt}$. 

This shows that the expression 
$(\sD,\lull)_{\pnt}:=(\del,\lull)_{\pnt}:\Oc_{\Curv,\pnt}^*\lra R^*$ 
is independent of the choice of local representative $\del$ at $q$. 

Let $h:R\lra S$ be a homomorphism of $k$-algebras. Then 
\begin{eqnarray*} 
\big(\Divfc{\Curv}(h)(\sD),\lull\big)_{\pnt} 
 & = & \pn{\Lin_h \circ \del,\lull}_{\pnt} \\ 
 & = & \Lin_h \circ (\del,\lull)_{\pnt} \\ 
 & = & \Lin_h \circ (\sD,\lull)_{\pnt} 
\end{eqnarray*} 
which implies \quad 
$ \fml \circ \, \Divfc{\Curv}(h) \,=\, \FDivf_{\Curv}(h) \circ \fml 
$.  
\ePf 

\bDef 
\label{Div_Z/C^0} 
Let $\Crvv$ be a projective curve over $k$,
and let $\pi:\Curv\lra\Crvv$ be its normalization. 
Then define the formal subgroup $\Divf_{\Curv/\Crvv}$ of $\Divfc{\Curv}$ 
to be the kernel of the composition $\pi_* \circ \fml$ 
\[ \Divf_{\Curv/\Crvv} = 
   \ker \Big( \Divfc{\Curv}
               \overset{\fml}\lra \FDivf_{\Curv} 
               \overset{\pi_*}\lra \FDivf_{\Crvv} 
        \Big) \laurin 
\] 
\eDef 

\bDef 
\label{Div_Y/X^0} 
Let $X$ be a projective variety over $k$, 
and $\pi: Y \lra X$ a birational morphism. 
%resolution of singularities. 
Then define the formal subgroup $\Divf_{Y/X}$ of $\Divfc{Y}$ as follows: 
\[  \Divf_{Y/X} =
    \bigcap_{C} \lrpn{ \llul\cut\Ct }^{-1} \Divf_{\Ct/C} 
\] 
where $C$ ranges over all Cartier curves in $X$ relative to $X \setminus U$ 
%the singular locus $X_{\sing}$ 
($\see$Def.~\ref{Cartier-curve}), 
where $\nu_C:\Ct\lra C$ is the normalization of $C$ 
and $\lul\cut \Ct:\Decfc{Y,\Ct} \lra \Divfc{\Ct}$ 
is the pull-back of relative Cartier divisors from $Y$ to $\Ct$ 
($\see$Propositions \ref{Decf formal group} and \ref{_.V}). 
Note that $\Ct \lra X$ factors through $Y$: 
let $C^Y$ denote the proper transform of $C$, 
i.e.\ the closure in $Y$ of $\pi^{-1}(C \isec U)$. 
As $\pi|_{C^Y}:C^Y \lra C$ is a birational morphism, 
the normalization $\Cyt$ of $\Cy$ is isomorphic to $\Ct$ 
and hence $\nu: \Ct \lra C$ 
factors through a morphism $\mu: \Ct \lra C^Y$. 
\eDef

\bRmk 
\label{coincidence old Div_Y/X}
For $\chr(k) = 0$, 
Definition \ref{Div_Z/C^0} of $\Divf_{\Curv/\Crvv}$ 
coincides with the definition of $\Divf_{\Curv/\Crvv}$ 
given in \cite[Proposition~3.23]{Ru}. 
Likewise, Definition \ref{Div_Y/X^0} of $\Divf_{Y/X}$ 
coincides with the definition of $\Divf_{Y/X}$ 
given in \cite[Proposition~3.24]{Ru}. 
\eRmk 

\bPf 
The second statement follows from the first by definition. 
For the first statement, observe that we used in Definition \ref{Div_Z/C^0} 
the fact that a relative Cartier divisor $\sD \in \Divfc{Z}\pn{R}$ 
induces a rational map $\del: Z \dra \Lin_R$ locally around each $\pnt \in Z(k)$ 
and we associated with $\sD$ the local symbols $\pn{\del,\lull}_{\pnt}$. 
For $\chr(k) = 0$, we are reduced to consider the cases 
$\Lin_R = \Gm$ and $\Lin_R = \Ga$. 
Then the explicit descriptions of the local symbol for $\Gm$ resp.\ $\Ga$ 
($\seecite$\cite[Example 3.13 resp.\ 3.14]{Ru}) 
yield the coincidence of the two definitions. 
\ePf

%\newpage 

\section{Category of Rtl.\ Maps Factoring through Rtl.\ Equivalence} 
\label{subsec:Mr_CH_0(X)_deg0}

We keep the notation fixed at the beginning of this note 
on page \pageref{prerequisites}: 
$X$ is a projective variety over $k$ 
which admits  a resolution of singularities $\pi:Y\lra X$, 
and $U\subset X_{\reg}$ a dense open subset, 
which we identify with its inverse image in $Y$. 
Algebraic groups are assumed to be smooth, connected and commutative, 
unless stated otherwise. 

%\newpage 

\bDef 
\label{MrCH(X)} 
Let $\MrCH{X}$ be the category of rational maps 
from $X$ to algebraic groups defined as follows: 
The objects of $\MrCH{X}$ are morphisms $\phe:U\lra G$ 
whose associated map on 0-cycles of degree 0 
\begin{eqnarray*}
  \Z_0(U)^0 & \lra & G(k) \\ 
  \sum l_i \, p_i & \lmt & \sum l_i \, \phe(p_i) 
\end{eqnarray*} 
factors through a homomorphism of groups $\CHSOo{X}\lra G(k)$. 
\footnote{A category of rational maps to algebraic groups is defined 
already by its objects, according to \cite[Remark 2.8]{Ru2}. } \\ 
%Remark~\ref{EquCatMr}. } \\ 
We refer to objects of $\MrCH{X}$ as 
rational maps from $X$ to algebraic groups 
\emph{factoring through rational equivalence} or 
\emph{factoring through $\CHSOo{X}$}. 
\eDef 

\bThm  
\label{Equivalence} 
The category $\MrCH{X}$ 
%of those rational maps from $X$ to algebraic groups 
%which factor through rational equivalence %$\CHSOo{X}$
is equal to the category $\Mr_{\Divor{Y/X}}$, 
i.e.\ to the category of those rational maps from $Y$ to algebraic groups 
which induce a transformation to $\Divor{Y/X}$ 
($\see$Definitions~\ref{Div_Z/C^0}, \ref{Div_Y/X^0}). 
\eThm 

\bPf  
A rational map from $Y$ to an algebraic group 
induces a transformation to $\Divor{Y/X}$ 
if and only if it induces a transformation to $\Divf_{Y/X}$ 
by \cite[Lemma 2.6]{Ru2}. %Proposition~\ref{induced Trafo}. 
Such a rational map is necessarily regular on $U$, 
since all $\sD \in \Divf_{Y/X}$ have support only on $Y\setminus U$. 
This follows from the fact that $\pi:Y\lra X$ is an isomorphism on $U$, 
and $\Divf_{Y/X}$ is the kernel of the push-forward $\pi_*$. 
Then according to Definitions \ref{CH_0(X)_deg0} 
and \cite[Definition 2.10]{Ru2} %\ref{Mr_sW} 
the task is to show that 
for a morphism $\phe:U\lra G$ from $U$ to an algebraic group $G$
with canonical decomposition $0\ra L\ra G\ra A\ra0$
the following conditions are equivalent: 

\begin{tabular}{rl}
(i) & $\phe\lrpn{\dv\pn{f}_{C}} = 0  \hspace{7.4em}
       \forall\lrpn{C,f}\in\fR_0\pn{X,U}$\laurink \\
(ii)& $\trafo_{\phe}\pn{\lam} \in \Divf_{Y/X} \pn{\Rfin} 
        \hspace{6em} \forall\,\Rfin\in\Artk,\; \lam\in\Ld(\Rfin)$\laurink 
\end{tabular} 

\noindent 
where $\trafo_{\phe}: \Ld \lra \Divfc{Y}$ 
is the induced transformation from \cite[Definition 2.4]{Ru2}. 
%Definition~\ref{induced Trafo}. 
Now if $h_{\lam}: G \lra \lam_*G$ denotes the push-out via 
$\lam \in \Ld(\Rfin) = \Homabk(L,\Lin_{\Rfin})$ and 
$\iota := \id_{\Lin_R} \in \Homabk(\Lin_{\Rfin},\Lin_{\Rfin}) 
 = \Lin_{\Rfin}^{\vee}(\Rfin)$, 
we obtain 
$\trafo_{\phe}\pn{\lam} = \trafo_{h_{\lam}\circ\phe}\pn{\iota} 
=: \dv_{\Rfin}\pn{h_\lam\circ\phe}$. 
This was explained in the proof of \cite[Theorem 2.12]{Ru2}. 
%\ref{Exist univObj}. 
By \cite[Lemma 1.12]{Ru2} %\ref{L subset L_S} 
obviously (i) is equivalent to 

\begin{tabular}{rll}
(i') & $\pn{h_{\lam} \circ \phe}\lrpn{\dv\pn{f}_{C}} = 0$  \hspace{4em} & 
       $\forall\lrpn{C,f}\in\fR_0\pn{X,U}$\laurink \\ 
     &  & $\forall \, \Rfin \in \Artk,\; \lam\in\Ld(\Rfin)$\laurin    
\end{tabular} 

\noindent 
Then the equivalence of (i) and (ii) follows from 
Lemma~\ref{L-bdl_factors_CH_0} below, 
if we replace $\phe$ in Lemma~\ref{L-bdl_factors_CH_0} 
by $h_{\lam} \circ \phe$. 
\ePf 

%\vspace{\vs} 

\bLem 
\label{L-bdl_factors_CH_0} 
Let $\phe:U\lra G$ be a morphism from $U$ 
to a smooth connected algebraic group $G\in\Ext\lrpn{A,\Lin_{\Rfin}}$, 
where $A$ is an abelian variety 
and $\Lin_{\Rfin}$ %an affine group. 
the linear group associated with $\Rfin\in\Artk$ 
($\seecite$\cite[Definition 1.11]{Ru2}). %($\see$Definition~\ref{Lin_R}). 
Then the following conditions are equivalent: 

\begin{tabular}{rl}
(i) & $\phe\lrpn{\dv\pn{f}_{C}} = 0    \hspace{8em} 
       \forall\lrpn{C,f}\in\fR_0\pn{X,U}$\laurink \\ 
(ii)& $\dv_{\Rfin}\pn{\phe}\in\Divf_{Y/X}\pn{\Rfin}$ \\
\end{tabular} 

\noindent 
where \;$\dv_{\Rfin}\pn{\phe} := \trafo_{\phe}\pn{\iota}$\; 
for $\trafo_{\phe}: \Lin_{\Rfin}^{\vee} \lra \Divfc{Y}$ 
the induced transformation from \cite[Definition 2.4]{Ru2} 
%Definition~\ref{induced Trafo} 
and $\iota := \id_{\Lin_R} \in \Homabk(\Lin_{\Rfin},\Lin_{\Rfin}) 
 = \Lin_{\Rfin}^{\vee}(\Rfin)$. 
\eLem 

%\newpage 

\bPf  (i)$\Llra$(ii) Let $C$ be a Cartier curve in
$X$ relative to $X\setminus U$, and let $\nu:\Ct\lra C$
be its normalization. 
Then Lemma~\ref{reciprocity} below asserts that 
the following conditions are equivalent:

\begin{tabular}{rl}
(j) & $\phe|_{C}\lrpn{\dv\pn{f}} = 0   \qquad\qquad\qquad
       \forall f\in\K\lrpn{C,U}^{*}$\laurink \\
(jj)& $\nu_{*}\lrpn{\dv_{\Rfin}\lrpn{\phe|_{C}}} = 0$\laurin 
\end{tabular} 

%\newpage

\noindent 
We have 
$\dv_{\Rfin}\lrpn{\phe|_{C}} = \dv_{\Rfin}\lrpn{\phe} \cut \Ct$, 
where $\lul\cut \Ct:\Decfc{Y,\Ct} \lra \Divfc{\Ct}$ 
is the pull-back of Cartier divisors from $Y$ to $\Ct$ 
($\see$Propositions \ref{Decf formal group} and \ref{_.V}). 
Then condition (i) is equivalent to 

\begin{tabular}{rl}
(iii) & $\lrpn{\nu_C}_{*}
        \big(\dv_{\Rfin}\lrpn{\phe} \cut \Ct \big) = 0  \hspace{18pt}
        \forall\,\textrm{Cartier curves }C \textrm{ relative to }
        X\setminus U$\laurin \\
\end{tabular} 

\noindent 
Conditions (iii) and (ii) are equivalent by definition of $\Divf_{Y/X}$
($\see$Def.s~\ref{Div_Z/C^0}, \ref{Div_Y/X^0}). 
%taking into account that 
%\;$\dv_{\Rfin}\lrpn{\phe}\in\Divorc{Y}\pn{\Rfin}$\; 
%according to \cite[Lemma 2.6]{Ru2}. %Proposition~\ref{induced Trafo}. 
\ePf 

%\newpage 

\bLem 
\label{reciprocity} 
Let $C$ be a projective curve and $\nu:\Ct\lra C$ its normalization. 
Let $\rmpC:C\dra G$ be a rational map from $C$ to a smooth connected 
algebraic group $G\in\Ext\lrpn{A,\Lin_{\Rfin}}$, 
i.e.\ $G$ is an $\Lin_{\Rfin}$-bundle over an abelian variety $A$, 
where $\Rfin$ is a finite $k$-algebra. 
Suppose that $\rmpC$ is regular on a dense open subset $U_{C}\subset C_{\reg}$,
which we identify with its preimage in $\Ct$. 
Then the following conditions are equivalent: 

\begin{tabular}{rl}
(i)  & $\rmpC\lrpn{\dv\pn{f}} = 0  
     \hspace{8.2em} \forall f\in\K\lrpn{C,U_{C}}^{*}$\laurink \\ 
(ii) & $\sum_{q\ra p} \big(\rmpC,g\big)_q = 0  
     \hspace{7em} \forall g\in\Oc_{C,p}^*\quad\forall p\in C(k)$\laurink \\
(iii) & $\lrpn{\nu_{*}\circ\fml} \lrpn{\dv_{\Rfin}\lrpn{\rmpC}} = 0$\laurin \\
\end{tabular} 
\eLem 

\bPf  (i)$\Llra$(ii) 
Let $f\in\K\lrpn{C,U_{C}}^{*}$. 
Write $\wt{f} := \nu^{\fis}f = f\circ\nu$.
Set $S:=\Ct\setminus U_{C}$. 
Using Lemma~\ref{normlzd} and the defining properties of a local
symbol ($\seecite$\cite[III, No.~1, Definition 2]{S}) 
%from Definition~\ref{locSymbol} 
we obtain 
\begin{eqnarray*}
\rmpC\,\big(\dv\pn{f}\big) 
 & = & \lrpn{\rmpC\circ\nu}\lrpn{\dv\big(\wt{f}\big)}\\
 & = & \sum_{c\notin S}\val_{c}\big(\wt{f}\big)\;\rmpC(c)\\
 & = & \sum_{c\notin S}\lrpn{\rmpC,\wt{f}}_{c}\\
 & = & -\sum_{s\in S}\lrpn{\rmpC,\wt{f}}_{s} \laurin 
% & = & -\sum_{s\in S}\lrpn{\left[\rmpC\right]_{\Phi_{s}},\wt{f}}_{s}\\
% & = & -\sum_{s\in S}\lrpn{\dv_{\Rfin}\lrpn{\rmpC},\wt{f}}_{s} 
\end{eqnarray*}
Let $\mdl$ be a modulus for $\rmpC$. 
For each $p\in\nu(S)$, each $g\in\Oc_{C,p}^*$ 
there is a rational function $f_{p}\in\K\lrpn{C,U_C}^{*}$
such that $\wt{g}\big/\wt{f}_{p}\big. \equiv 1 \mod \mdl$ at $\nu^{-1}(p)$ 
and $\wt{f}_{p}\equiv 1\mod\mdl$ at $s$ for all $s\in S\setminus\nu^{-1}(p)$, 
by the approximation theorem. 
Indeed, for $\mdl = \sum_{s\in S} n_s s$  and given $g\in\Oc_{C,p}^*$ 
the approximation theorem yields an element $b_p \in \sK_{\Ct}$ 
such that $b_p \equiv \wt{g} \mod \mdl$ at $\nu^{-1}(p)$ 
and $b_p \equiv 1\mod \mdl$ at $s$ for all $s \in S\setminus\nu^{-1}(p)$. 
Then $b_p$ can be chosen such that $b_p = \wt{f}_p$ %$b_p = \nu^{\fis} f_p$ 
for some $f_p \in \K\lrpn{C,U_C}^{*}$. 
This means that $g = f_p + h = f_p \lrpn{1 + f_p^{-1} h}$ 
for some $h \in \mc_{C,p}$ with $\wt{h} \in \mc_{\Ct,q}^{n_q}$ for each $q \ra p$. 
As $f_p \in \sO_{C,p}^*$, we have $\wt{f}_p^{-1} \wt{h} \in \mc_{\Ct,q}^{n_q}$, 
which implies the statement. 

Then $\big(\rmpC,\wt{f}_p\big)_q = \big(\rmpC,\wt{g}\big)_q$ 
for all $q\in\nu^{-1}(p)$ and 
$\big(\rmpC,\wt{f}_p\big) = 0$ for all $s\in S\setminus\nu^{-1}(p)$, 
by \cite[III, No.~1, Definition 2, (i),(ii)]{S}. 
%by properties (a) and (c) of Definition~\ref{locSymbol}. 
Hence $\rmpC\big(\dv\pn{f_{p}}\big) = 0$ if and only if 
$\sum_{q\ra p}\big(\rmpC,\wt{f}_{p}\big)_q = 
 \sum_{q\ra p}\big(\rmpC,\wt{g}\big)_q = 0$. 
Thus $\rmpC\big(\dv\pn{f}_{C}\big) = 0$
for all $f\in\K\lrpn{C,U_{C}}^{*}$ if and only if 
$\sum_{q\ra p}\big(\rmpC,\wt{g}\big)_q = 0$
for all $g\in\Oc_{C,p}^*$, $p\in\nu(S)$. 
It remains to remark that $\big(\rmpC,h\big)_c = 0$ for all 
$h\in\Oc_{\Ct,c}^*\supset\Oc_{C,\nu(c)}^*$,
$c\in U_{C}$, since $\rmpC$ and $h$ are both regular at $c$. 

(ii)$\Llra$(iii) Let $p\in C$. 
For each $q\in\nu^{-1}(p)$ let $\Phi_{q}:U_{q}\tms L \lra G_Y$
be a local trivialization of the induced $L$-bundle over $\Ct$
in a neighbourhood $U_{q}\ni q$. 
For all $f\in\sO_{\Ct,q}^*$ we have 
\begin{eqnarray*} 
\big(\rmpC,f\big)_q & = & \big(\left[\rmpC\right]_{\Phi_{q}},f\big)_q \\ 
                   & = & \big(\dv_{\Rfin}\lrpn{\rmpC},f\big)_q 
\end{eqnarray*} 
by \cite[Lemma 3.16]{Ru} %Lemma~\ref{locSymb_triv} 
and as was seen in the proof of Proposition~\ref{fml}. 
Then %the condition 
$\sum_{q\ra p}\lrpn{\rmpC,\wt{g}}_q = 0$ for all $g\in\Oc_{C,p}^*$ 
is equivalent to the condition that the image 
$\sum_{q\ra p}\big(\dv_{\Rfin}\lrpn{\rmpC},\lull\big)_q$ of 
$\dv_{\Rfin}\lrpn{\rmpC}$ in 
$\dsum_{q\ra p}\Homabk\lrpn{\Oc_{\Ct,q}^*,\Rfin^*}$ 
vanishes on $\Oc_{C,p}^*$, %by construction, 
which says 
\[ 0 = \sum_{q\ra p}\big(\dv_{\Rfin}\lrpn{\rmpC},\lull\big)_q 
       \circ\Bigpn{\wh{\nu_q^{\fis}}}^* 
       \in\dsum_{q\ra p}\Homabk\lrpn{\Oc_{C,\nu(q)}^*,\Rfin^*} \laurin 
\] 
This is true for all $p\in C$ if and only if 
\;$\lrpn{\nu_{*}\circ\fml} \big(\dv_{\Rfin}\lrpn{\rmpC}\big) = 0$\; 
by definition of the push-forward for formal divisors 
($\see$Definition~\ref{phi_*}). 
\ePf

\section{Kernel of the Push-forward of Divisors} 
\label{sub:Div_Y/X^0}

The functor $\Divf_{Y/X}$ was introduced in Section~\ref{FormalDivisors} 
as ``the kernel of the push-forward of relative divisors'' 
($\see$Definitions \ref{Div_Z/C^0}, \ref{Div_Y/X^0}). 
The goal of this section is to show that 
\;$\Divor{Y/X} := \Divf_{Y/X} \tms_{\Picf_Y} \Picorf{Y}$\; 
is represented by a dual-algebraic formal group. 

\bNot 
For any subfunctor $\Fc$ of $\Divf_Y$ 
we denote by $F^{0,\red}$ the induced subfunctor of $\cl^{-1}\Picorf{Y}$: 
\[  \Fc^{0,\red} = \Fc \tms_{\Picf_Y} \Picorf{Y} \laurin 
\] 
\eNot 

\bRmk 
\label{restr fml grp} 
%If $\fmlG$ is a subfunctor of $\Divfc{Y}$ which is a formal group, 
%then we denote by $\fmlG^{0,\red}$ the induced subfunctor of $\cl^{-1}\Picorf{Y}$: 
%\[ \fmlG^{0,\red} = \fmlG \tms_{\wh{\Picf_Y}} \wh{\Picorf{Y}} \laurin 
%\] 
If $\fmlG$ is a formal group, then $\fmlG^{0,\red}$ is again a formal group, 
since $\fmlG^{0,\red} = \fmlG \tms_{\wh{\Picf_Y}} \wh{\Picorf{Y}}$ 
and the category $\Gfk$ of formal groups admits fibre-products 
($\seecite$\cite[I, 4.9, p.~35]{Fo}). 
\eRmk 

%\newpage 

\bLem 
\label{fml inj} 
Let $\Curv$ be a normal curve over $k$. 
Then the natural transformation \;$\fml:\Divfc{\Curv} \lra \FDivf_{\Curv}$\; 
induces a monomorphism of formal $k$-groups 
\;$\Divorc{\Curv} \lra \FDivf_{\Curv}$. 
\eLem 

\bPf 
Let $R$ be a finite $k$-algebra and $\sD\in\Divorc{\Curv}(R)$ with 
$\fml\lrpn{\sD}=0$. 
%The divisor $\sD$ arises from the canonical local section 
%of the $R^*$-bundle $\Pb_{R^*}(\sD)$. 
Let $\Gp_{\Lin_R}(\sD)\in\Extabk(\Alba{\Curv},\Lin_R)
                  \cong\Pic^0_{\Alba{\Curv}}(R)$ 
($\seecite$\cite[Proposition 1.19]{Ru2}) %($\see$Proposition \ref{Xeil}) 
be the image of $\sO_{\Curv \tens R}(\sD)$ under the homomorphism 
\[ \alb: \Picor{\Curv} \lra \Alb\big(\Picor{\Curv}\big) 
                       = \Pic^0\lrpn{\Pic^0\big(\Picor{\Curv}\big)} 
                       = \Pic^0 \bigpn{\Alba{\Curv}} \laurin 
\] 
Then $\Pb_{\Lin_R}(\sD) := \Gp_{\Lin_R}(\sD)\tms_{\Alba{\Curv}}\Curv$ 
is the $\Lin_R$-bundle on $\Curv$ 
associated with $\sO_{\Curv \tens R}(\sD)$. 
The canonical 1-section of $\sO_{\Curv \tens R}(\sD)$ 
induces a rational map 
\[ \rmp{\sD}:\Curv \dra \Pb_{\Lin_R}(\sD)
            = \Gp_{\Lin_R}(\sD)\tms_{\Alba{\Curv}}\Curv \lra \Gp_{\Lin_R}(\sD) 
\] 
%\;$ \rmp{\sD}: \Curv \dra \Gp_{\Lin_R}(\sD) 
%$\; 
s.t.\ $\sD = \trafo_{\rmp{\sD}}\pn{\iota}$, where 
$\trafo_{\rmp{\sD}}: \Lin_R^{\vee} \lra \Divfc{\Curv}$ 
is the induced transformation from \cite[Definition 2.4]{Ru2} 
%Definition~\ref{induced Trafo} 
and 
$\iota := \id_{\Lin_R} \in \Homabk(\Lin_R,\Lin_R) = 
  \Lin_R^{\vee}(R)$. 
If $g\in\Oc_{\Curv,\pnt}^*$ for some  $\pnt \in \Curv(k)$, 
then by \cite[Lemma~3.16]{Ru} 
%Lemma \ref{locSymb_triv} 
the local symbol $\pn{\rmp{\sD},g}_q$ lies in the fibre of $\Gp_{\Lin_R}$ 
over $0 \in \Alba{\Curv}$, which is $\Lin_R$, 
and $\pn{\rmp{\sD},g}_q = \pn{\sD,g}_q$. 
As $\fml\pn{\sD}=0$, it holds 
$0 = \pn{\sD,g}_q = \pn{\rmp{\sD},g}_q$ 
for all $\pnt\in\Curv(k)$ and all $g\in\Oc_{\Curv,\pnt}^*$.  
Lemma~\ref{reciprocity} implies that 
$0 = \rmp{\sD}\bigpn{\dv\pn{f}}$ 
for all $f\in\K\lrpn{\Curv,\Curv_{\reg}}^*$, 
i.e.\ $\rmp{\sD}$ factors through $\CHSOo{\Curv}$. 
But $\CHSOo{\Curv} = \Alba{\Curv}$ is an abelian variety, 
since $Z_{\sing} = \varnothing$, 
thus $\rmp{\sD}$ extends to a morphism defined on all of $\Curv$. 
The transformation $\trafo_{\rmp{\sD}}$ has the property that 
$\Supp\pn{\im \trafo_{\rmp{\sD}}} \subset Z$ is the locus 
where $\rmp{\sD}$ is not defined. 
Then $\im \pn{\trafo_{\rmp{\sD}}} = 0$, %i.e. $\trafo_{\rmp{\sD}}$ is the zero-map, 
hence $\sD = \trafo_{\rmp{\sD}}\pn{\iota} = 0$. 

\hfill 
\ePf 

\bThm 
\label{Div_Z/C repr}
Let $\Crvv$ be a projective curve over $k$, 
and let $\pi:\Curv\lra\Crvv$ be its normalization. 
Then the functor $\Divor{\Curv/\Crvv}$ is represented 
by a dual-algebraic formal group. 
\eThm 

\bPf 
First note that the normalization $\pi:\Curv\lra\Crvv$ is an isomorphism 
on the preimage $\pi^{-1}\Crvv_{\reg}$ of the regular locus of $\Crvv$. 
Hence the push forward $\pi_*$ is an isomorphism on formal divisors 
supported on $\pi^{-1}\Crvv_{\reg}$. 
Therefore the formal divisors in \;$\ker\big(\FDivf_{\Curv}\lra\FDivf_{\Crvv}\big)$ 
have support only on the inverse image $\Sing_{\Curv}$ 
of the singular locus $\Sing$ of $\Crvv$, which is finite. 
\begin{eqnarray*}
& & \ker\Big(\FDivf_{\Curv}\lra\FDivf_{\Crvv}\Big) \\ 
&=& \ker\lrpn{ \dsum_{\pnt\in\Curv(k)} \Homfabk\lrpn{\Oc_{\Curv,\pnt}^*,k^*} 
    \lra \dsum_{\pntt\in\Crvv(k)} \Homfabk\lrpn{\Oc_{\Crvv,\pntt}^*,k^*} }\\ 
&=& \dsum_{\pntt\in\Crvv(k)}
    \ker\lrpn{ \dsum_{\pnt\ra\pntt} \Homfabk\lrpn{\Oc_{\Curv,\pnt}^*,k^*} 
               \lra \Homfabk\lrpn{\Oc_{\Crvv,\pntt}^*,k^*} }\\ 
&=& \dsum_{\pntt\in\Crvv(k)}
    \ker\bigg( \Homfabk\lrpn{\Oc_{\Curv,\pntt}^*,k^*} 
               \lra \Homfabk\lrpn{\Oc_{\Crvv,\pntt}^*,k^*} \bigg)\\ 
&=& \dsum_{\pntt\in\Sing(k)} 
    \Homfabk\Big(\left.\Oc_{\Curv,\pntt}^*\right/\Oc_{\Crvv,\pntt}^*\;,\;k^*\Big)\\ 
&=& \dsum_{\pntt\in\Sing(k)} 
    \Homfabk\Big(\left.\sO_{\Curv,\pntt}^*\right/\sO_{\Crvv,\pntt}^*\;,\;k^*\Big)\\ 
&=& \Homfabk\Big(\left.\sO_{\Curv}^*\right/\sO_{\Crvv}^*\;,\;k^*\Big) 
\end{eqnarray*} 
Now $\left.\sO_{\Curv}^*\right/\sO_{\Crvv}^* = 
     \prod_{\pntt\in\Sing}\left.\sO_{\Curv,\pntt}^*\right/\sO_{\Crvv,\pntt}^*$ 
is a sheaf of abelian groups over $\Algk$ 
with Lie-algebra 
\;$\Lie \lrpn{\left. \sO_{\Curv}^* \right/ \sO_{\Crvv}^*} 
 = \left. \sO_{\Curv} \right/ \sO_{\Crvv} 
 = \prod_{\pntt\in\Sing} \left. \sO_{\Curv,\pntt} \right/ \sO_{\Crvv,\pntt}$. 
%Considering this object as a coherent sheaf on a projective curve 
This is a coherent sheaf 
concentrated on finitely many points ${\pntt\in\Sing}$, 
which implies that $\Lie\lrpn{\left.\sO_{\Curv}^*\right/\sO_{\Crvv}^*}$ 
is finite dimensional. 
Then the group sheaf $\left.\sO_{\Curv}^*\right/\sO_{\Crvv}^*$ 
is represented by an affine algebraic group 
%$L \cong \prod_{\pntt\in\Sing} \Lin_{R_{\pntt}}$ 
$L \cong \prod_{\pntt\in\Sing} \Trz_{\pntt} \tms \Upf_{R_{\pntt}}$, 
%where $\Lin_{R_{\pntt}}$ is the linear group 
%from \cite[1.1.6]{Ru2} %\ref{LinGroup_Ring} 
%associated with the finite $\fld$-algebra 
%$R_{\pntt} = \left. \sO_{\Curv,\pntt} \right/ \sO_{\Crvv,\pntt}$. 
where $\Trz_{\pntt}$ denotes the torus $\bigpn{\prod_{\pnt \ra \pntt} \Gm} / \Gm$
%where $\GmS{\pnt}$ denotes a group isomorphic to the multiplicative group 
%and attached to the point $\pnt \in \Curv(\fld)$, 
and $\Upf_{R_{\pntt}}$ is the unipotent group 
from \cite[1.1.6]{Ru2} %\ref{LinGroup_Ring} 
associated with the finite $\fld$-algebra 
$R_{\pntt} = \fld + \bigpn{\dsum_{\pnt \ra \pntt} \fm_{\Curv,\pnt}} / \fm_{\Crvv,\pntt}$. 
Then $\ker\lrpn{\FDivf_{\Curv}\lra\FDivf_{\Crvv}} = 
\Homfabk\lrpn{\left.\sO_{\Curv}^*\right/\sO_{\Crvv}^*,k^*} = \Ld$ 
is the Cartier-dual of an affine algebraic group, 
hence a dual-algebraic formal group. 

$\Divfc{\Curv}$ is a formal group 
by \cite[Proposition~2.1]{Ru2}, %Proposition~\ref{Divf formal group}, 
hence $\Divorc{\Curv}$ is again a formal group, 
according to Remark~\ref{restr fml grp}. 
Due to Lemma \ref{fml inj} the transformation 
\;$\fml^{0,\red}:\Divorc{\Curv} \lra \FDivf_{\Curv}$\; 
is a monomorphism of formal groups. 
%and so is the homomorphism $\Ld\lra\FDivf_{\Curv}$. 
Then the fibre-product of formal groups 
\;$\Divor{\Curv/\Crvv} = \Divorc{\Curv} \tms_{\FDivf_{\Curv}} \Ld$\; 
is a subgroup of the dual-algebraic formal group $\Ld$, 
hence dual-algebraic by \cite[Lemma~3.15]{Ru2}. 
%($\see$Lemma~{subFmlGroup_algebraic}). 
\ePf 

\bThm 
\label{Div_Y/X repr} 
Let $X$ be a projective variety over $k$ 
that admits a resolution of singularities $\pi: Y \lra X$. 
Then the functor $\Divor{Y/X}$ is represented 
by a dual-algebraic formal group. 
\eThm 

The proof of Theorem~\ref{Div_Y/X repr} 
is similar to the proof of \cite[Proposition 3.24]{Ru}. 
But instead of considering $k$-valued points and Lie-algebra of the functor 
$\Divor{Y/X}$ separately (which makes sense only for $\chr(k) = 0$), 
we can control the values of this functor for all finite rings $\Rfin$ 
by a formal group $\Fm{Y}{\mdl}$ ($\seecite$\cite[Definition 3.13]{Ru2}) 
%($\see$Definition~\ref{DefFm(X,D)}) 
associated with a sufficiently large effective divisor 
$\mdl$ on $Y$ with support on the inverse image 
$\SgY$ of the singular locus $\Sing=X_{\sing}$ of $X$. 
This theorem will be fundamental 
for an independent proof of the existence of $\Urq{X}$. 
Conversely, 
Theorem \ref{Div_Y/X repr} will follow from the existence of $\Urq{X}$ 
(proven in \cite[Theorem 1]{ESV}) and Proposition \ref{Urq-Divoyx}. 
%Therefore we leave it at that declarations. 
\medskip 

\bPf [Proof of Thm~\ref{Div_Y/X repr}] 
The functor 
$\Divor{Y/X} = 
\bigcap_{C} \Decfc{Y,\Ct} \tms_{\Divfc{\Ct}} \Divor{\Ct/C}$, 
as a projective limit of formal groups, is again a formal group. 
(A formal group functor is a formal group if and only if 
it commutes with finite projective limits, and projective limits commute.) 
It remains to show that $\Divor{Y/X}$ is dual-algebraic. 
According to \cite[Proposition 3.21]{Ru2} %Proposition~\ref{exhaust Div_Y^0} 
it is sufficient to show that there exists an effective divisor $\mdl$ on $Y$ 
such that $\Divor{Y/X} \subset \Fm{Y}{\mdl}$. 

Let $\sD \in \Divf_Y\pn{\Rfin}$ 
be a non-trivial relative Cartier divisor on $Y$ 
for some $\Rfin \in \Artk$. 
Assume $\Supp\pn{\sD}$ is not contained in the inverse image 
$\SgY=\Sing\tms_X Y$ of the locus $\Sing = X \setminus U$. 
Then $\pi\big(\Supp(\sD)\big)$ on $X$ is not contained in $\Sing$. 
Let $\sL$ be a very ample line bundle on $X$, 
consider the space $|\sL|^{d-1}$, where $d = \dim X$, 
of complete intersection curves $C = H_1 \cap \ldots \cap H_{d-1}$ 
with $H_i \in |\sL| = \Prj\big(\H^0(X,\sL)\big)$ 
for $i = 1,\ldots,d-1$. 
For Cartier curves $C$ in $|\sL|^{d-1}$ 
%relative to $X_{\sing}$ 
(with $C_Y = C \tms_X Y$) 
the following properties are open and dense: 

\begin{tabular}{rl} 
(a) & $C$ intersects $\pi\bigpn{\Supp(\sD)} \cap U$ properly\laurin \\ 
(b) & $\sD \cut C_Y$ is a non-trivial divisor on $C_Y \cap \pi^{-1}U$\laurin \\ 
(c) & $C_Y$ is regular in a neighbourhood of 
$\Supp\pn{\sD \cut C_Y} \cap \pi^{-1}U$\laurin 
\end{tabular} 

\noindent 
(a) is due to the fact that $\sL$ is very ample, 
(b) follows from (a) and the fact that $\Supp(\sD)$ 
is locally a prime divisor at almost every $q\in Y$ 
and (c) is a consequence of the Bertini theorems. 
Therefore there exists a Cartier curve $C$ in $X$ 
satisfying the conditions (a)--(c). 
As the normalization $\nu:\Ct\lra C$ coincides with $\pi|_{C_Y}$, 
it is an isomorphism on $C_Y \cap U$. 
Thus $\nu_*\bigpn{\sD\cut\Ct} \neq 0$. 
This implies $\sD \notin \Divor{Y/X}\pn{\Rfin}$. 
Hence $\Supp\bigpn{\Divor{Y/X}}$ is contained in $\SgY$. 

Let $\Es$ be the reduced effective divisor associated with the sum of 
the components of $\SgY$ of codimension 1 in $Y$. 
If $C$ is a Cartier curve in $X$ relative to $\Sing$, 
then the normalization $\nu: \Ct \lra C$ 
factors through the proper transform $\Cy$. 
%the image $C_Y$ of $\mu$ consists of those components 
%of the inverse image $C \tms_X Y$ that are not contained 
%in the exceptional divisor of the desingularization $Y \lra X$. 
%Assume $\chr(k) = p > 0$. 
%If $\sL$ is a sufficiently ample, %line bundle on $X$, 
%a sufficiently high power of an ample line bundle on $X$, 
By Lemma~\ref{cover property} 
one finds a family $T \subset |\sL|^{d-1}$ of Cartier curves $C$ 
such that the set $\bigcup_{C \in T} C^Y \cap \Es$ 
contains an open dense subset $W$ of $\Es$. 
%This is due to the fact that $Y$ is obtained from $X$ 
%by repeatedly blowing up. 
Since $\Divor{\Ct/C}$ is dual-algebraic, 
by \cite[Proposition 3.21]{Ru2} %Proposition~\ref{exhaust Div_Y^0} 
there exists a natural number $n_C \in \Nat$ such that 
$\Divor{\Ct/C} \subset \Fm{\Ct}{n_C \Es\cut\Ct}$. 
As $\Divor{\Ct/C}$ is Cartier dual to the affine part of $\Pic^0_{C}$, 
the bound $n_c$ 
is controlled by the dimension of the unipotent part of $\Pic^0_{C}$. 
(It is sufficient that $n_C$ satisfies $\fm_{\Ct,q}^{n_C} \subset \fm_{C,p}$ 
for all $p\in C$, $q\in\nu^{-1}(p)$, 
see \cite[Section 9.2, proof of Proposition 9]{BLR}). 
By upper semi-continuity of $\dim \Pic^0_{C}$ 
for the curves $C$ in $|\sL|^{d-1}$, 
we may assume that there exists a common bound 
$n \geq n_C$ for all $C \in T$. 
If $\chr(k)=p>0$ 
and $p$ divides $n$, replace $n$ by $n+1$. 
Now $\sD \in \Divor{Y/X}\pn{\Rfin}$ 
implies by definition 
$\sD\cut\Ct \in \Divor{\Ct/C}\pn{\Rfin}$ for all Cartier curves $C$. 
Then $\sD\cut\Ct \in \Fm{\Ct}{n\Es\cut\Ct}$ for all $C \in T$. 
By construction of $T \subset |\sL|^{d-1}$ and $n$, 
this yields $\sD \in \Fm{Y}{n\Es}$. 
\ePf 

\bLem 
\label{cover property}
Let $X$ be a projective variety over $k$ 
which admits a projective resolution of singularities $\pi: Y \lra X$. 
%Then the following holds: 
%\[ \bigcup_{C} \, \Cy \,\ni\, \pnt  \hspace{20mm} 
%   \forall \, \pnt \in Y  \textrm{ of codim 1} 
%\]  
Then \;$\bigcup_{C} \, \Cy = Y$,\, 
i.e.\ the union of the proper transforms of the curves $C$ in $X$ 
covers $Y$, %up to codimension 2, 
where $C$ ranges over all Cartier curves 
in $X$ relative to $X \setminus U$. 
%
%\begin{tabular}{rl} 
%$(\dagger)$ & $\bigcup_{C} \Cy \ni \pnt  \hspace{10mm} 
%\forall \, \pnt \in Y \textrm{ of codim 1}$\laurink \\ 
% & where $C$ ranges over all Cartier curves 
%     in $X$ relative to $X \setminus U$\laurin 
%\end{tabular}
\eLem 

\bPf 
According to \cite[II, Theorem~7.17]{H} we may assume that 
$\pi:Y \lra X$ is given as the blowing up of an ideal sheaf on $X$, 
%i.e.\ $Y = \Proj \dsum_{\nu \geq 0} \sJ^{\nu}$ 
%i.e.\ $Y$ is a relative projective space over $X$ and 
in particular $\pi_* \sO_Y(1) = \sJ$ is an ideal sheaf on $X$. 
Let $\sL$ be an ample line bundle on $X$. 
Then $\Ampl :=  \sO_Y(1) \tens \, \pi^*\sL^N$ 
for sufficiently 
large integer $N$ 
is a very ample line bundle on $Y$ over $X$ 
($\seecite$\cite[II, Proposition 7.10 (b)]{H})
and hence  over $k$, since $X$ is projective. 

For a given closed point $\pnt \in Y$ 
we are now going to construct a Cartier curve $C$ on $X$ 
such that $\pnt \in \Cy$. 
Choose a tangent vector $\tg \in \Tg_{\pnt} Y \setminus \Tg_{\pnt} S_Y$, 
where $\Tg_{\pnt} Y$ resp.\ $\Tg_{\pnt} S_Y$ 
denotes the tangent space of $Y$ resp.\ $S_Y$ at $\pnt$, 
%that is the fibre of the tangent sheaf $\Tgs_V$ of $V$ at $\pnt$, 
%for any variety $V$. 
and $S_Y$ is the inverse image in $Y$ of the locus of $S = X \setminus U$. 
%This is possible since $\dim \Tg_{\pnt} Y = \dim \Oma_Y \tens k(\pnt) = \dim Y$ 
%for all $\pnt \in Y$ by smoothness of $Y$. 
Since $\Ampl$ is very ample, 
one finds a complete intersection curve $C_{\pnt,\tg} \in |\Ampl|^{d-1}$, 
where $d = \dim X$, 
through $\pnt$ such that $\tg \in \Tg_{\pnt} C_{\pnt,\tg}$. 
The curve $C_{\pnt,\tg} = \bigcap_{i=1}^{d-1} \Zero\pn{s_i}$ 
is given as intersection of the divisors of zeroes $\Zero\pn{s_i}$ for 
\begin{eqnarray*}
s_i & \in & \Gam\bigpn{Y, \Ampl} \\ 
 & = & \Gam\bigpn{X, \pi_*\Ampl} \\ 
 & = & \Gam\Bigpn{X, \pi_*\bigpn{\pi^*\sL^N \tens_{\sO_Y} \sO_Y(1)}} \\ 
 & = & \Gam\bigpn{X, \sL^N \tens_{\sO_X} \pi_*\sO_Y(1)} \\ 
 & = & \Gam\bigpn{X, \sL^N \tens_{\sO_X} \sJ} \\ 
 & \subset & \Gam\bigpn{X, \sL^N} \laurin 
\end{eqnarray*}
%where we used projection formula in the third equation. 
%and the fact that $\sJ$ is an ideal sheaf on $X$ in the last line. 
Thus $C_{\pnt,\tg}$ is the proper transform of the curve 
$C = \bigcap_{i=1}^{d-1} \Zero\pn{s_i} \in |\sL^N|^{d-1}$ on $X$, 
which is a Cartier curve by construction. 
\ePf

%\newpage 

\section{Universal Reg.\ Quotient of the Chow Group} 
\label{subsec:Universal_Regular_Quotient} 

The results obtained up to now provide the necessary foundations for a 
functorial description of the generalized Albanese variety $\Urq{X}$ 
of Esnault-Srinivas-Viehweg 
%for a (singular) projective variety $X$
and its dual, 
which was the initial intention of this work.

\subsection{Existence and Construction} 
\label{subsubsec:Exist_UnivRegQuot} 

%The universal regular quotient $\Urq{X}$ of $\CHSOo{X}$ 
The Albanese $\Urq{X}$ of Esnault-Srinivas-Viehweg 
and the rational map \;$\urq{X}: X \dra \Urq{X}$\; 
are by definition the universal object for
the category $\MrCH{X}$ of morphisms from $U \subset X_\reg$ 
factoring through $\CHSOo{X}$ ($\see$Definition~\ref{MrCH(X)}). 
%which factor through a homomorphism of groups $\CHSOo{X}\lra G$. 
Therefore $\Urq{X}$ is often called 
\emph{the universal regular quotient of $\CHSOo{X}$}. 

\bPrp 
\label{Urq-Divoyx}
The following conditions are equivalent: 

\begin{tabular}{rl} 
(i) & $\Urq{X}$ exists as an algebraic group\laurin \\ 
(ii) & $\Divor{Y/X}$ is represented by a dual-algebraic formal group\laurin 
\end{tabular} 
\ePrp 

\bPf 
The proof of \cite[Theorem 2.12]{Ru2} %Theorem~\ref{Exist univObj} 
shows that a category $\Mr$ of rational maps from $Y$ to algebraic groups 
satisfying \cite[$2.3.1, \lrpn{\diamondsuit \;1,2}$]{Ru2}
admits a universal object if and only if the formal group $\fmlG \subset \Divfc{Y}$ 
generated by $\bigcup_{\phe \in \Mr} \im \pn{\trafo_{\phe}}$ is dual-algebraic 
and $\Mr = \Mr_{\fmlG}$. 
By Theorem~\ref{Equivalence} the category $\MrCH{X}$ is equal 
to the category $\Mr_{\Divor{Y/X}}$ of rational maps 
from a desingularization $Y$ of $X$ to algebraic groups 
which induce a transformation to the formal group $\Divor{Y/X}$. 
Thus $\MrCH{X}$ admits a universal object if and only if 
$\Divor{Y/X}$ is dual-algebraic. 
\ePf 

\vspace{\vs} 

\bPf[Proof of Theorem~\ref{main_result}] 
The existence of $\Urq{X}$ was proven in \cite[Theorem 1]{ESV}. 
Alternatively we showed the existence %of $\Urq{X}$ 
from a direct proof of Theorem~\ref{Div_Y/X repr} 
in combination with Proposition~\ref{Urq-Divoyx}. 

According to Theorem~\ref{Equivalence} it holds 
$\Urq{X} = \Alb_{\Divor{Y/X}}(Y)$, 
and the latter %$\Alb_{\Divor{Y/X}}(Y)$ 
was constructed as the dual 1-motive of $\bigbt{\Divor{Y/X} \lra \Picorf{Y}}$ 
($\seecite$\cite[Remark 2.14]{Ru2}). %($\see$Remark~\ref{Alb_constr}). 
This gives the existence and an explicit construction 
of the universal regular quotient of $\CHSOo{X}$, 
as well as a description of its dual. 
The proof of Theorem~\ref{main_result} 
is thus complete. 
\ePf

%\newpage 

\subsection{Case without Desingularization} 
\label{no Desingularization}

If the characteristic of the base field $k$ is positive, 
it is not known whether there exists always a desingularization 
of the given projective variety $X$. 
Therefore we want to get rid of this assumption, 
%the resolution of singularities of $X$ 
which was used 
in the functorial description of the universal regular quotient of $\CHSOo{X}$. 
In this subsection we show that we can replace the desingularization of $X$ 
by any birational projective morphism $Y \ra X$ 
where the variety $Y$ is normal and has the property that 
every rational map from $Y$ to an abelian variety is regular on all of $Y$. 
Then we construct a variety $Y$ with the required properties by blowing up 
an ideal sheaf on $X$. 

\bThm 
\label{dualUrq}
Let $\pi: Y \lra X$ be a projective morphism with the properties 

\begin{tabular}{rl} 
$\lrpn{\bir \, 1}$ & $\pi: Y \lra X$ is birational\laurin \\ 
$\lrpn{\bir \, 2}$ & Every rational map from $X$ to an abelian variety $A$ \\ 
 & extends to a morphism from $Y$ to $A$. \\ 
 $\lrpn{\bir \, 3}$ & $\Picor{Y}$ is an abelian variety. 
\end{tabular} 

\noindent 
Then the functor $\bigbt{\Divor{Y/X} \lra \Picorf{Y}}$ 
is dual (in the sense of 1-motives) to the universal regular quotient 
$\Urq{X}$ of $\CHSOo{X}$. 
\eThm 

\bPf 
Let $\Fctav{X}$ denote the category of rational maps from $X$ to abelian varieties, 
$\Morav{X}$ the category of morphisms from $X$ to abelian varieties. 
Let $\Albfct{X}$ denote the classical Albanese variety of $X$ 
in the sense of functions, 
i.e.\ the universal object for $\Fctav{X}$ ($\seecite$\cite{La}), 
and $\Albmor{X}$ the Albanese variety of $X$ in the sense of morphisms, 
i.e.\ the universal object for $\Morav{X}$ ($\seecite$\cite{S2}). 
Conditions $\lrpn{\bir \, 1}$ and $\lrpn{\bir \, 2}$ imply that 
$\Fctav{X}$ is equivalent to $\Morav{Y}$, 
hence $\Albfct{X} = \Albmor{Y}$. 
Since $\Picor{Y}$ is an abelian variety by condition $\lrpn{\bir \, 3}$, 
the functoriality of $\Picor{}$ yields a universal mapping property 
of $\Pic^0 \bigpn{\Picor{Y}}$ for $\Morav{Y}$, 
hence $\Albmor{Y}$ is the dual of $\Picor{Y}$. 

Now the concept of 
\emph{categories of rational maps from $Y$ to algebraic groups} 
can be set up in the same way as it was done for smooth proper varieties in 
%Subsection \ref{sub:Categories-of-Rational} 
\cite[Subsection~2.2]{Ru2}: 
A rational map $\phe: Y \dra G$ from $Y$ to an algebraic group $G$ 
with affine part $L$ and abelian quotient $A$ 
induces a transformation $\trafo_{\phe}: \Ld \lra \Divfc{Y}$ as in 
%Subsubsection \ref{sec:indTrafo} 
\cite[Subsubsection~2.2.1]{Ru2}, 
here we only needed $A\lrpn{\sK_{Y}} = A\lrpn{\sO_{Y,\pnt}}$ 
for every point $\pnt \in Y$ 
(which is nothing else than to say that every rational map from $Y$ to $A$ 
is regular at every $\pnt \in Y$, 
and this is satisfied by conditions $\lrpn{\bir \, 1}$ and $\lrpn{\bir \, 2}$). 
Then the existence and construction of universal objects 
for categories of rational maps from $Y$ to algebraic groups as in 
%Subsubsection \ref{subsub:Exist+Construct} 
\cite[Subsubsection~2.2.2]{Ru2}, 
the equality of the categories $\MrCH{X}$ and $\Mr_{\Divor{Y/X}}$ as in 
Theorem~\ref{Equivalence} 
and the representability of $\Divor{Y/X}$ as in Theorem~\ref{Div_Y/X repr} 
carry over literally. 
Hence the proof of Theorem~\ref{main_result} from above is valid 
in this situation as well. 
\ePf 

\bThm 
\label{extendRatlMap}
To a given projective variety $X$ over $k$ 
there exists always a projective morphism $\pi: Y \lra X$ 
satisfying the properties $\lrpn{\bir \, 1 - 3}$ from Theorem~\ref{dualUrq}. 
\eThm 

\bPf 
Consider the universal rational map $\alp: X \dra \Albfct{X}$\
from $X$ to the Albanese of $X$ in the sense of functions. 
As $\Albfct{X}$ is an abelian variety, it is projective 
($\seecite$\cite[II, Theorem~8.12]{P}). 
Let $\emb: \Albfct{X} \lra \Prj^r$ be an embedding of $\Albfct{X}$ 
into projective space. 
Then the composition $\emb \circ \alp: X \dra \Prj^r$ 
lifts to a morphism $\bet: V \lra \Prj^r$, 
where $\sig: V \lra X$ is the blowing up of a suitable ideal sheaf $\sJ$ on $X$ 
($\seecite$\cite[II, Example~7.17.3]{H}). 
The closed subscheme corresponding to $\sJ$ has support equal to 
$X \setminus U$, 
where $U$ is the maximal open set on which $\alp$ 
%\;$\alb: X \dra \Albfct{X}$\; 
is regular. 
\[ \xymatrix{ V \ar[d]_{\sig} \ar[r]^{\bet}  &  \Prj^r \\ 
                     X \ar@{-->}[r]^-{\alp}  &  \Albfct{X} \ar[u]_{\emb} 
   }
\] 
The image of $\bet$ is the closure of the image of $\iota \circ \alp$ in $\Prj^r$, 
which is contained in the image of $\Albfct{X}$ in $\Prj^r$: 
\[ \bet\pn{V} \,\subset\, \ol{\bet\pn{\sig^{-1}U}} \,=\, \ol{\emb \alp\pn{U}} 
   \,\subset\, \emb\bigpn{\Albfct{X}} \laurin 
\] 
Thus $\bet$ induces a morphism $\alp_V: V \lra \Albfct{X}$ 
extending $\alp \circ \sig|_{\sig^{-1}U}$. 
Since $\alp: X \dra \Albfct{X}$ is defined by the universal mapping property 
for rational maps from $X$ to abelian varieties, 
every rational map $\phe: X \dra A$ factors as 
$\phe = h \circ \alp: X \dra \Albfct{X} \lra A$ 
and hence extends to a morphism $\phe_V = h \circ \alp_V: V \lra \Albfct{X} \lra A$. 
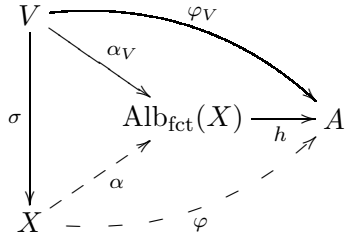
\[ \xymatrix{ V \ar[dd]_{\sig} \ar[dr]^{\alp_V} \ar@/^1.3pc/[drr]^{\phe_V}  &  & \\ 
    & \Albfct{X} \ar[r]_-{h \;}  &  A \\ 
   X \ar@{-->}[ur]_{\alp} \ar@{-->}@/_1.3pc/[urr]_{\phe}  &  & \\ 
   }
\] 
Let $\nu: Y \lra V$ be the normalization of $V$. 
A morphism from $V$ to $A$ yields a morphism from $Y$ to $A$ 
via pull-back $\nu^*$. 
Then $\pi:= \sig \circ \nu: Y \lra V \lra X$ satisfies conditions $\lrpn{\bir \, 1,2}$. 
According to \cite[Theorem~5.4]{K}, 
the identity component $\Pic_{Y}^0$ of the Picard scheme is projective. 
As the base field $k$ is perfect, 
the underlying reduced scheme $\Picor{Y}$ is an abelian variety, 
which gives condition $\lrpn{\bir \, 3}$. 
\ePf 

\vspace{\vs} 

\bPf[Proof of Theorem~\ref{main_result_2}]
Follows directly from Theorem \ref{dualUrq} and Theorem \ref{extendRatlMap}. 
\ePf

%\newpage 

\subsection{Descent of the Base Field} 
\label{UnivRegQuot_descent}

Let $\fld$ be a perfect field. 
%Let $\unidom$ be a universal domain, 
%i.e.\ $\unidom$ is an algebraically closed extension of $\fld$ 
%of infinite transcendence degree over $\fld$. 
Let $\clfld$ be an algebraic closure of $\fld$. %(in $\unidom$). 
Let $\Xo$ be a projective variety defined over $\fld$. 
%let $Y \lra X$ be a projective resolution of singularities over $\fld$ 
%(resp.\ an alteration over $\fld$ such that $Y$ is regular and projective). 
Write $\Xb = X \tens_{\fld} \clfld$. %and $\Yb = Y \tens_{\fld} \clfld$. 

We are going to show that the field of definition 
of the universal regular quotient $\Urq{\Xb}$ of $\CHSOo{\Xb}$ 
descends to $\fld$. 
When we do not assume that $\Xo$ is endowed with a $\fld$-rational point, 
the Albanese of Esnault-Srinivas-Viehweg for $\Xo$ 
exists only as a $\fld$-torsor $\Urqq{1}{\Xo}$ 
for an algebraic $\fld$-group $\Urqq{0}{\Xo}$. 

\bDef 
\label{difference map}
If $\phe: \Xo \dra P$ is a rational map to a torsor $P$ for an algebraic group $G$
we define $\phe^{(-)}: \Xo \tms \Xo \dra G$ to be 
the rational map which assigns to $(p,q) \in \Xo \tms \Xo$ 
the unique $g \in G$ such that $g \cdot \phe(p) = \phe(q)$. 
We set \;$\phe^{(1)} := \phe$\; and \;$\phe^{(0)} := \phe^{(-)}$. 
%(cf.\ \cite[Notation 2.18]{Ru2}). 
%($\see$Notation~\ref{phe^(i)}). 
\eDef 

\bDef 
\label{factor geom through CHo} 
%Suppose $\anyfld$ is an extension of finite type of $\fld$ in $\unidom$ 
%and let $\clanyfld$ be the algebraic closure (in $\unidom$). 
Let $\phe: \Xo \dra P$ be a rational map 
defined over $\fld$ to a $\fld$-torsor $P$ 
for a $\fld$-group $G$. 
Then %the base-changed  map 
$\phe\tens_{\fld}\clfld: \Xo\tens_{\fld}\clfld 
 \dra P\tens_{\fld}\clfld$ 
over $\clfld$ maps 0-cycles of degree 0 to $G\pn{\clfld}$: 
%If $z\in\Z_0(\Xo\tens_{\fld}\anyfld)_{\deg 0}$ is a 0-cycle of degree 0, 
%we may consider $\phe\lrpn{z}$ as an element of $G(\anyfld)$: 
%Over the albegraic closure $\clanyfld$ the cycle 
a cycle $z \in \Z_0(\Xo\tens_{\fld}\clfld)^0$ 
decomposes into a sum of 0-cycles of the form $q-p$. 
By the property of $P\tens_{\fld}\clfld$ of being a torsor for 
$G\tens_{\fld}\clfld$, one may consider 
$\pn{\phe\tens_{\fld}\clfld} \lrpn{z}$ as an element of 
$G\pn{\clfld}$. %(cf.\ Notation~\ref{Difference map}). 

We say that $\phe$ \emph{factors geometrically through rational equivalence}, 
if the base-changed map $\phe\tens_{\fld}\clfld$ 
%$\phe\tens_{\anyfld}\clanyfld: \Xb\tens_{\clfld}\clanyfld 
% \dra P\tens_{\anyfld}\clanyfld$ 
over $\clfld$ factors through rational equivalence in the sense of 
Definition~\ref{MrCH(X)}. 
\eDef 

%\newpage 

\bThm 
\label{descent_UnivRegQuot} 
There exists a $\fld$-torsor $\Urqq{1}{\Xo}$ for an algebraic 
$\fld$-group $\Urqq{0}{\Xo}$ and rational maps defined over $\fld$ 
\[ \urqqo{i}: \Xo^{2-i} \dra \Urqq{i}{\Xo} 
\] 
for $i=1,0$, satisfying the following universal property: 

%Suppose $\anyfld$ is an extension of finite type of $\fld$ in $\unidom$ 
%and let $\clanyfld$ be the algebraic closure (in $\unidom$). 
If $\phe: \Xo \dra \Trr{1}$ 
is a rational map defined over $\fld$ to a $\fld$-torsor $\Trr{1}$ 
for an algebraic $\fld$-group $\Trr{0}$, 
factoring geometrically through rational equivalence 
%$\CHSOo{\Xo\tens_{\fld}\anyfld}$ 
($\see$Definition~\ref{factor geom through CHo}), 
there exist a unique affine homomorphism of $\fld$-torsors $\homm{1}$ 
and a unique homomorphism of algebraic $\fld$-groups $\homm{0}$ 
\[ \homm{i}: \Urqq{i}{\Xo} \lra \Trr{i} 
\] 
defined over $\fld$, such that 
\[ \phe^{(i)} = \homm{i} \circ \urqqo{i} 
\] 
for $i=1,0$. 
Here $\phe^{(i)}$ are the rational maps from 
Definition~\ref{difference map}. 
%Here $\phe^{(1)} = \phe$ and $\phe^{(0)}: \Xo \tms \Xo \lra G^{(0)}$ 
%is the rational map which assigns to $(p,q) \in \Xo \tms \Xo$ 
%the unique $g \in G^{(0)}$ such that it holds 
%$g \cdot \phe(p) = \phe(q)$ 
%(cf.\ \cite[Notation 2.21]{Ru2}). 
%($\see$Notation~\ref{phe^(i)}). 
\eThm 

\bPf 
Since $\Xo$ and the singular locus $\Xo_{\sing}$ are defined over $k$, 
the absolute Galois group $\Gal\pn{\clfld/\fld}$ acts on $\Xb$ and on 
the group $\CHSOo{\Xb}$. 
In particular, for the conjugates by means of $\gal \in \Gal\pn{\clfld/\fld}$
one has the equalities 
$\Xb^{\gal} = \Xb$ and $\bigpn{\CHSOo{\Xb}}^{\gal} = \CHSOo{\Xb}$.  

Let $\pn{\urq{\Xb}}^{\gal} : \Xb \dra \Urq{\Xb}^{\gal}$ 
be the transform of \;$\urq{\Xb}: \Xb \dra \Urq{\Xb}$ 
by means of $\gal \in \Gal\pn{\clfld/\fld}$. 
Since $\pn{\urq{\Xb}}^{\gal}$ factors through 
$\bigpn{\CHSOo{\Xb}}^{\gal} = \CHSOo{\Xb}$, 
it is an object of $\MrCH{\Xb}\pn{\Xb}$. 
Then the map $\pn{\urq{\Xb}}^{\gal}$ factors as 
\;$\pn{\urq{\Xb}}^{\gal} = \hg{\gal} \circ \urq{\Xb}$, 
where \;$\hg{\gal}: \Urq{\Xb} \lra \Urq{\Xb}^{\gal}$\; is an affine homomorphism 
(i.e.\ a homomorphism of algebraic groups composed with a translation). 
Now we can apply the machinery of Galois descent, 
which is described e.g.\ in \cite[V, \S~4]{S}, cf.\ \cite[2.3.3]{Ru2}. 
%\ref{descentBaseField}. 
In particular, the proof of Theorem~\ref{descent_UnivRegQuot} 
is the same as in \cite[V, No.~22]{S}. 
\ePf

\newpage 
%\quad 
%\newpage

\section*{Glossary} 
\label{Glossar}
%\addcontentsline{toc}{section}{Glossary} 

\hspace{+0.01em} 
{\large Categories} 
\medskip \\  
\begin{tabular}{lll} 
$\Set$ & sets & \cite[1.1.1]{Ru2} \\ %\ref{k-Func Compl} \\ 
$\Ab$ & abelian groups & \cite[1.1.1]{Ru2} \\ %\ref{k-Func Compl} \\ 
$\Algk$ & $\fld$-algebras & \cite[1.1.1]{Ru2} \\ %\ref{k-Func Compl} \\ 
$\Artk$ & $\fld$-algebras of finite length & \cite[1.1.1]{Ru2} \\ %\ref{k-Func Compl} \\ 
$\Fctr(\mathsf{A},\mathsf{B})$ & functors from $\mathsf{A}$ to $\mathsf{B}$ & 
\cite[1.1.1]{Ru2} \\ %\ref{k-Func Compl} \\ 
$\Abkp$ & $\fld$-group functors ($= \Fctr(\Algk,\Ab)$) & \cite[1.1.1]{Ru2} \\ 
%\ref{k-groupFunc}\\ 
$\Abk$ & $\fld$-group sheaves (for fppf-topology) & \cite[1.1.5]{Ru2} \\ 
%\ref{Sheaves}\\ 
$\Gk$ & $\fld$-groups (= $\fld$-group schemes) & \cite[1.1.3]{Ru2} \\ 
%\ref{sub:Algebraic-Groups} \\ 
$\Gak$ & affine $\fld$-groups & \cite[1.1.3]{Ru2} \\ %\ref{sub:Algebraic-Groups} \\ 
$\aGk$ & algebraic $\fld$-groups & \cite[1.1.3]{Ru2} \\ 
%\ref{sub:Algebraic-Groups} \\ 
$\aGak$ & affine algebraic $\fld$-groups & \cite[1.1.3]{Ru2} \\ 
%\ref{sub:Algebraic-Groups} \\ 
$\Gfk$ & formal $\fld$-groups & \cite[1.1.4]{Ru2} \\ %\ref{sub:Formal-Groups} \\ 
%$\fGfk$ & formal $\fld$-groups of finite type & \cite[1.1.4]{Ru2} \\ 
%\ref{sub:Formal-Groups} \\ 
$\dGfk$ & dual-algebraic formal $\fld$-groups & \cite[1.2.1]{Ru2} \\ 
%\ref{Anatomy} \\ 
%$\Compl(\Abk)$ & complexes of $\fld$-group functors & 1.4 \\ 
%\ref{1-MotiveUnip} \\ 
%$\Mot$ & 1-motives $M = (\fmlG,L,A,G,\mu)$ & 1.4.3 \\ %\ref{Category} \\ 
%$\Mk$ & $[\fmlG\ra G]$ with $\fmlG\in\dGfk$, $G\in\aGk$ & 1.4.3 \\ 
%\ref{Category} \\ 
%$\Mkqi$ & localization of $\Mk$ at quasi-isomorphisms & 1.4.3 \\ 
%\ref{Category} \\ 
\end{tabular} 

%\newpage 

\vspace{\vsvs} 
\noindent 
\hspace{+0.01em} 
{\large Functors} 
\medskip\\ 
\begin{tabular}{lll} 
$\wh{\Fc}$ & completion of $\Fc \in \Abkp$ & \cite[1.1.1]{Ru2} \\ 
%\ref{k-Func Compl} \\ 
$\fmlG_{\et}$ & $= \fmlG \circ \red$ \;\'etale part of $\fmlG \in \Gfk$ & 
\cite[1.1.4]{Ru2} \\ 
%\ref{k-groupFunc} \\ 
$\fmlG_{\inf}$ & $= \ker\bigpn{\fmlG \ra \fmlG_{\et}}$ 
                         infinitesimal part of $\fmlG \in \Gfk$ & \cite[1.1.4]{Ru2} \\ 
                         %\ref{k-groupFunc} \\ 
$\Lin_R$ & $= \Gm\pn{\llul\tens R}$ linear group ass.\ to $R\in\Algk$ & 
\cite[1.1.6]{Ru2} \\ %\ref{sub:Lin_R} \\ 
$\Trz_R$ & $= \Lin_{R_{\red}}$ torus ass.\ to $R\in\Algk$ & \cite[1.1.6]{Ru2} \\ 
%\ref{sub:Lin_R} \\ 
$\Upf_R$ & $= \ker\bigpn{\Lin_R \ra \Trz_R}$ 
                            unipotent group ass.\ to $R\in\Algk$ & \cite[1.1.6]{Ru2} \\ 
                            %\ref{sub:Lin_R} \\ 
$\Fm{X}{\mdl}$ & $\subset \Divf_X$ formal group of modulus $\mdl$ & 
\cite[3.2.1]{Ru2} \\ % \ref{Alb(X,D)exist} 
$\Picf_X$ & Picard functor of $X$ & \cite[2.1]{Ru2} \\ %\ref{Pic-functor} \\ 
$\Picorf{X}$ & $=$ reduced identity component of $\Picf_X$ & \cite[2.1]{Ru2} \\ 
%\ref{Pic-functor} \\ 

$\Divf_X$ & relative Cartier divisors on $X$ & \ref{FuncDiv} \\ 
$\Divor{X}$ & $= \cl^{-1}\Picorf{X}$ & \ref{sub:Div_Y/X^0} \\ 
%$\Divf_X^S$ & $\subset \Divf_X$ with support in $S$ & \ref{FuncDiv} \\ 
$\Decf_{X,V}$ & $\subset \Divf_X$ intersecting $V$ properly & \ref{FuncDiv} \\ 
$\FDivf_C$ & formal divisors on a curve $C$ & \ref{FormalDivisors} \\ 
$\Divf_{Z/C}$ & $= \ker\bigpn{\Divfc{Z} \ra \FDivf_Z \ra \FDivf_C}$ & 
\ref{FormalDivisors} \\ %\ref{sub:Div_Y/X^0} \\ 
$\Divf_{Y/X}$ & $= \bigcap_C \bigpn{\llul\cut\Ct}^{-1} \Divf_{\Ct/C}$ 
\; if $Y \ra X$ is birational & \ref{FormalDivisors} \\ %\ref{sub:Div_Y/X^0} \\ 
%$\Divszc$ & $= \ker \Bigpn{ \Divs{Z} \lra \FDivf_Z \lra \FDivf_C }$ & 
%\ref{no Desingularization} \\ 
$\Divsyx$ & $= \bigcap_C \lrpn{ \llul \cut \Cyt }^{-1} \Divf_{\Ct/C} \cut \, \Cyt$ & 
%\ref{no Desingularization} \\ 
\end{tabular} 

\newpage 
%\vspace{\vsvs} 
\noindent 
\hspace{+0.01em} 
{\large Algebraic groups} 
\medskip\\ 
\begin{tabular}{lll} 
$\Jacc{\Crv}$ & Jacobian of a curve $C$ & \cite[3.3]{Ru2} \\ %\ref{subsec:JacMod} \\ 
$\Jacm{\Crv}{\mdl}$ & Jacobian of $C$ of modulus $\mdl$ & \cite[3.3]{Ru2} \\ 
%\ref{subsec:JacMod} \\ 
$\Lm{\Crv}{\mdl}$ & affine part of $\Jacm{C}{\mdl}$ & \cite[3.3]{Ru2} \\ 
%\ref{subsec:JacMod} \\ 
$\Tm{\Crv}{\mdl}$ & torus part of $\Jacm{C}{\mdl}$ & \cite[3.3]{Ru2} \\ 
%\ref{subsec:JacMod} \\ 
$\Um{\Crv}{\mdl}$ & unipotent part of $\Jacm{C}{\mdl}$  & \cite[3.3]{Ru2} \\ 
%\ref{subsec:JacMod} \\ 
$\Pic_X$ & Picard scheme of $X$ & \cite[2.1.2]{Ru2} \\ %\ref{Pic-functor} \\ 
$\Picor{X}$ & Picard variety of $X$ & \cite[2.1.2]{Ru2} \\ %\ref{Pic-functor} \\ 
 & ($=$ reduced identity component of $\Pic_X$) & \\ 
$\Alba{X}$ & Albanese variety of $X$ %, if $X$ smooth 
 & \cite[2.3.1]{Ru2} \\ %\ref{subsub:Exist+Construct} \\ 
$\AlbF{X}$ & $= \lrbt{\fmlG \ra \Pic_X^0}^{\vee}$ 
                       universal object for $\Mr_{\fmlG}$ & \cite[2.3.1]{Ru2} \\ 
                       %\ref{subsub:Exist+Construct} \\ 
$\AlbbF{1}{X}$ & universal torsor for $\Mr_{\fmlG}$ & \cite[2.3.3]{Ru2} \\ 
%\ref{descentBaseField} \\
$\AlbbF{0}{X}$ & universal group for $\Mr_{\fmlG}$ & \cite[2.3.3]{Ru2} \\ 
%\ref{descentBaseField} \\
$\Albm{X}{\mdl}$ & Albanese variety of $X$ of modulus $\mdl$ & \cite[3.2.1]{Ru2} \\ 
%\ref{Alb(X,D)exist} \\ 
$\Urq{X}$ & universal regular quotient of $\CHSOo{X}$ & 
\ref{subsubsec:Exist_UnivRegQuot} \\ 
%$\Albsg{X}$ & $= \Urq{X}$ & \ref{sec:Exm_Surf} \\ 
$\Urqq{1}{X}$ & torsor descended from $\Urq{\Xb}$ & 
\ref{UnivRegQuot_descent} \\ 
$\Urqq{0}{X}$ & group descended from $\Urq{\Xb}$ & 
\ref{UnivRegQuot_descent} \\ 
\end{tabular} 

\vspace{\vsvs} 
\noindent 
\hspace{+0.01em} 
{\large Chow Groups of 0-cycles} 
\medskip\\ 
\begin{tabular}{lll} 
$\CHm{X}{\mdl}$ & relative Chow group of $X$ of modulus $\mdl$ & 
\cite[3.4]{Ru2} \\ %\ref{subsec:ChowMod} \\ 
$\CHmo{X}{\mdl}$ & $= \ker\bigpn{\deg: \CHm{X}{\mdl} \ra \Zint}$ & 
\cite[3.4]{Ru2} \\ %\ref{subsec:ChowMod} \\ 
$\CHSO{X}$ & %$= \CHm{X}{X_{\sing}}$ 
relative Chow group by Levine-Weibel & \ref{Chow_LW} \\ 
$\CHSOo{X}$ & $= \ker\bigpn{\deg: \CHSO{X} \ra \Zint}$ & \ref{Chow_LW} \\ 
\end{tabular}

\vspace{\vsvs} 
\noindent 
\hspace{+0.01em} 
{\large Rational Maps} 
\medskip\\ 
\begin{tabular}{lll} 
$\Mr$ & a category of rational maps & \cite[2.2]{Ru2} \\ 
%\ref{sub:Categories-of-Rational} \\ 
$\MrCH{X}$ & rational maps factoring through $\CHSOo{X}$ &  \ref{subsec:Mr_CH_0(X)_deg0}%{MrCH(X)} 
\\ 
$\MrCHm{X}{\mdl}$ & rational maps factoring through $\CHmo{X}{\mdl}$ & 
\cite[3.4]{Ru2} \\ %\ref{subsec:ChowMod}%{MrCH(X,D)} \\ 
$\Mrm{X}{\mdl}$ & $= \bigst{\phe\,|\,\modu\pn{\phe} \leq \mdl}$ & 
\cite[3.2.1]{Ru2} \\ %\ref{Alb(X,D)exist} \\ 
$\Mr_{\fmlG}$ & $= \bigst{\phe\,|\,\im\pn{\trafo_{\phe}} \subset \fmlG}$ & 
\cite[2.2]{Ru2} \\ %\ref{sub:Categories-of-Rational} \\ 
$\trafo_{\phe}$ & $: \Ld \ra \Divfc{X}$ induced transformation of $\phe$ & 
\cite[2.2]{Ru2} \\ %\ref{sub:Categories-of-Rational} \\ 
$\modu\pn{\phe}$ & modulus of rational map $\phe$ & \cite[3.2]{Ru2} \\ 
%\ref{subsec:AlbMod} \\ 
%$$ &  & \ref{} \\ 
%$$ &  & \ref{} \\ 
\end{tabular}

\end{document}